\documentclass[12pt]{amsart}
\usepackage{amssymb}
\begin{document}
\theoremstyle{plain}
\newtheorem{thm}{Theorem}[section]
\newtheorem{theorem}[thm]{Lemma}
\newtheorem{lemma}[thm]{Lemma}
\newtheorem{corollary}[thm]{Corollary}
\newtheorem{proposition}[thm]{Proposition}
\theoremstyle{definition}
\newtheorem{notations}[thm]{Notations}
\newtheorem{remark}[thm]{Remark}
\newtheorem{remarks}[thm]{Remarks}
\newtheorem{definition}[thm]{Definition}
\newtheorem{claim}[thm]{Claim}
\newtheorem{assumption}[thm]{Assumption}
\numberwithin{equation}{thm}
\newcommand{\zar}{{\rm zar}}
\newcommand{\an}{{\rm an}}
\newcommand{\red}{{\rm red}}
\newcommand{\codim}{{\rm codim}}
\newcommand{\rank}{{\rm rank}}
\newcommand{\Pic}{{\rm Pic}}
\newcommand{\Div}{{\rm Div}}
\newcommand{\Hom}{{\rm Hom}}
\newcommand{\im}{{\rm Im}}
\newcommand{\Spec}{{\rm Spec}}
\newcommand{\sing}{{\rm sing}}
\newcommand{\reg}{{\rm reg}}
\newcommand{\Char}{{\rm char}}
\newcommand{\Tr}{{\rm Tr}}
\newcommand{\Gal}{{\rm Gal}}
\newcommand{\Min}{{\rm Min \ }}
\newcommand{\Max}{{\rm Max \ }}
\newcommand{\soplus}[1]{\stackrel{#1}{\oplus}}
\newcommand{\dlog}{{\rm dlog}\,}    
\newcommand{\limdir}[1]{{\displaystyle{\mathop{\rm
lim}_{\buildrel\longrightarrow\over{#1}}}}\,}
\newcommand{\liminv}[1]{{\displaystyle{\mathop{\rm
lim}_{\buildrel\longleftarrow\over{#1}}}}\,}
\newcommand{\boxtensor}{{\Box\kern-9.03pt\raise1.42pt\hbox{$\times$}}}
\newcommand{\sext}{\mbox{${\mathcal E}xt\,$}}
\newcommand{\shom}{\mbox{${\mathcal H}om\,$}}
\newcommand{\coker}{{\rm coker}\,}
\renewcommand{\iff}{\mbox{ $\Longleftrightarrow$ }}
\newcommand{\onto}[1]{\mbox{$\,\>#1 >>\hspace{-.5cm}\to\hspace{.15cm}$}}
\catcode`\@=11
\def\opn#1#2{\def#1{\mathop{\kern0pt\fam0#2}\nolimits}}
\def\bold#1{{\bf #1}}%
\def\underrightarrow{\mathpalette\underrightarrow@}
\def\underrightarrow@#1#2{\vtop{\ialign{$##$\cr
 \hfil#1#2\hfil\cr\noalign{\nointerlineskip}%
 #1{-}\mkern-6mu\cleaders\hbox{$#1\mkern-2mu{-}\mkern-2mu$}\hfill
 \mkern-6mu{\to}\cr}}}
\let\underarrow\underrightarrow
\def\underleftarrow{\mathpalette\underleftarrow@}
\def\underleftarrow@#1#2{\vtop{\ialign{$##$\cr
 \hfil#1#2\hfil\cr\noalign{\nointerlineskip}#1{\leftarrow}\mkern-6mu
 \cleaders\hbox{$#1\mkern-2mu{-}\mkern-2mu$}\hfill
 \mkern-6mu{-}\cr}}}
\let\amp@rs@nd@\relax
\newdimen\ex@
\ex@.2326ex
\newdimen\bigaw@
\newdimen\minaw@
\minaw@16.08739\ex@
\newdimen\minCDaw@
\minCDaw@2.5pc
\newif\ifCD@
\def\minCDarrowwidth#1{\minCDaw@#1}
\newenvironment{CD}{\@CD}{\@endCD}
\def\@CD{\def\A##1A##2A{\llap{$\vcenter{\hbox
 {$\scriptstyle##1$}}$}\Big\uparrow\rlap{$\vcenter{\hbox{%
$\scriptstyle##2$}}$}&&}%
\def\V##1V##2V{\llap{$\vcenter{\hbox
 {$\scriptstyle##1$}}$}\Big\downarrow\rlap{$\vcenter{\hbox{%
$\scriptstyle##2$}}$}&&}%
\def\={&\hskip.5em\mathrel
 {\vbox{\hrule width\minCDaw@\vskip3\ex@\hrule width
 \minCDaw@}}\hskip.5em&}%
\def\verteq{\Big\Vert&&}%
\def\noarr{&&}%
\def\vspace##1{\noalign{\vskip##1\relax}}\relax\let\amp@rs@nd@&\iffalse}\fi
 \CD@true\vcenter\bgroup\relax\let\\=\cr\iffalse}\fi\tabskip\z@skip\baselineskip20\ex@
 \lineskip3\ex@\lineskiplimit3\ex@\halign\bgroup
 &\hfill$\m@th##$\hfill\cr}
\def\@endCD{\cr\egroup\egroup}
\def\>#1>#2>{\amp@rs@nd@\setbox\z@\hbox{$\scriptstyle
 \;{#1}\;\;$}\setbox\@ne\hbox{$\scriptstyle\;{#2}\;\;$}\setbox\tw@
 \hbox{$#2$}\ifCD@
 \global\bigaw@\minCDaw@\else\global\bigaw@\minaw@\fi
 \ifdim\wd\z@>\bigaw@\global\bigaw@\wd\z@\fi
 \ifdim\wd\@ne>\bigaw@\global\bigaw@\wd\@ne\fi
 \ifCD@\hskip.5em\fi
 \ifdim\wd\tw@>\z@
 \mathrel{\mathop{\hbox to\bigaw@{\rightarrowfill}}\limits^{#1}_{#2}}\else
 \mathrel{\mathop{\hbox to\bigaw@{\rightarrowfill}}\limits^{#1}}\fi
 \ifCD@\hskip.5em\fi\amp@rs@nd@}
\def\<#1<#2<{\amp@rs@nd@\setbox\z@\hbox{$\scriptstyle
 \;\;{#1}\;$}\setbox\@ne\hbox{$\scriptstyle\;\;{#2}\;$}\setbox\tw@
 \hbox{$#2$}\ifCD@
 \global\bigaw@\minCDaw@\else\global\bigaw@\minaw@\fi
 \ifdim\wd\z@>\bigaw@\global\bigaw@\wd\z@\fi
 \ifdim\wd\@ne>\bigaw@\global\bigaw@\wd\@ne\fi
 \ifCD@\hskip.5em\fi
 \ifdim\wd\tw@>\z@
 \mathrel{\mathop{\hbox to\bigaw@{\leftarrowfill}}\limits^{#1}_{#2}}\else
 \mathrel{\mathop{\hbox to\bigaw@{\leftarrowfill}}\limits^{#1}}\fi
 \ifCD@\hskip.5em\fi\amp@rs@nd@}
\newenvironment{CDS}{\@CDS}{\@endCDS}
\def\@CDS{\def\A##1A##2A{\llap{$\vcenter{\hbox
 {$\scriptstyle##1$}}$}\Big\uparrow\rlap{$\vcenter{\hbox{%
$\scriptstyle##2$}}$}&}%
\def\V##1V##2V{\llap{$\vcenter{\hbox
 {$\scriptstyle##1$}}$}\Big\downarrow\rlap{$\vcenter{\hbox{%
$\scriptstyle##2$}}$}&}%
\def\={&\hskip.5em\mathrel
 {\vbox{\hrule width\minCDaw@\vskip3\ex@\hrule width
 \minCDaw@}}\hskip.5em&}
\def\verteq{\Big\Vert&}
\def\novarr{&}
\def\noharr{&&}
\def\SE##1E##2E{\slantedarrow(0,18)(4,-3){##1}{##2}&}
\def\SW##1W##2W{\slantedarrow(24,18)(-4,-3){##1}{##2}&}
\def\NE##1E##2E{\slantedarrow(0,0)(4,3){##1}{##2}&}
\def\NW##1W##2W{\slantedarrow(24,0)(-4,3){##1}{##2}&}
\def\slantedarrow(##1)(##2)##3##4{%
\thinlines\unitlength1pt\lower 6.5pt\hbox{\begin{picture}(24,18)%
\put(##1){\vector(##2){24}}%
\put(0,8){$\scriptstyle##3$}%
\put(20,8){$\scriptstyle##4$}%
\end{picture}}}
\def\vspace##1{\noalign{\vskip##1\relax}}\relax\let\amp@rs@nd@&\iffalse}\fi
 \CD@true\vcenter\bgroup\relax\let\\=\cr\iffalse}\fi\tabskip\z@skip\baselineskip20\ex@
 \lineskip3\ex@\lineskiplimit3\ex@\halign\bgroup
 &\hfill$\m@th##$\hfill\cr}
\def\@endCDS{\cr\egroup\egroup}
\newdimen\TriCDarrw@
\newif\ifTriV@
\newenvironment{TriCDV}{\@TriCDV}{\@endTriCD}
\newenvironment{TriCDA}{\@TriCDA}{\@endTriCD}
\def\@TriCDV{\TriV@true\def\TriCDpos@{6}\@TriCD}
\def\@TriCDA{\TriV@false\def\TriCDpos@{10}\@TriCD}
\def\@TriCD#1#2#3#4#5#6{%
\setbox0\hbox{$\ifTriV@#6\else#1\fi$}
\TriCDarrw@=\wd0 \advance\TriCDarrw@ 24pt
\advance\TriCDarrw@ -1em
\def\SE##1E##2E{\slantedarrow(0,18)(2,-3){##1}{##2}&}
\def\SW##1W##2W{\slantedarrow(12,18)(-2,-3){##1}{##2}&}
\def\NE##1E##2E{\slantedarrow(0,0)(2,3){##1}{##2}&}
\def\NW##1W##2W{\slantedarrow(12,0)(-2,3){##1}{##2}&}
\def\slantedarrow(##1)(##2)##3##4{\thinlines\unitlength1pt
\lower 6.5pt\hbox{\begin{picture}(12,18)%
\put(##1){\vector(##2){12}}%
\put(-4,\TriCDpos@){$\scriptstyle##3$}%
\put(12,\TriCDpos@){$\scriptstyle##4$}%
\end{picture}}}
\def\={\mathrel {\vbox{\hrule
   width\TriCDarrw@\vskip3\ex@\hrule width
   \TriCDarrw@}}}
\def\>##1>>{\setbox\z@\hbox{$\scriptstyle
 \;{##1}\;\;$}\global\bigaw@\TriCDarrw@
 \ifdim\wd\z@>\bigaw@\global\bigaw@\wd\z@\fi
 \hskip.5em
 \mathrel{\mathop{\hbox to \TriCDarrw@
{\rightarrowfill}}\limits^{##1}}
 \hskip.5em}
\def\<##1<<{\setbox\z@\hbox{$\scriptstyle
 \;{##1}\;\;$}\global\bigaw@\TriCDarrw@
 \ifdim\wd\z@>\bigaw@\global\bigaw@\wd\z@\fi
 \mathrel{\mathop{\hbox to\bigaw@{\leftarrowfill}}\limits^{##1}}
 }
 \CD@true\vcenter\bgroup\relax\let\\=\cr\iffalse}\fi
 \tabskip\z@skip\baselineskip20\ex@
 \lineskip3\ex@\lineskiplimit3\ex@
 \ifTriV@
 \halign\bgroup
 &\hfill$\m@th##$\hfill\cr
#1&\multispan3\hfill$#2$\hfill&#3\\
&#4&#5\\
&&#6\cr\egroup%
\else
 \halign\bgroup
 &\hfill$\m@th##$\hfill\cr
&&#1\\%
&#2&#3\\
#4&\multispan3\hfill$#5$\hfill&#6\cr\egroup
\fi}
\def\@endTriCD{\egroup}
\newcommand{\sA}{{\mathcal A}}
\newcommand{\sB}{{\mathcal B}}
\newcommand{\sC}{{\mathcal C}}
\newcommand{\sD}{{\mathcal D}}
\newcommand{\sE}{{\mathcal E}}
\newcommand{\sF}{{\mathcal F}}
\newcommand{\sG}{{\mathcal G}}
\newcommand{\sH}{{\mathcal H}}
\newcommand{\sI}{{\mathcal I}}
\newcommand{\sJ}{{\mathcal J}}
\newcommand{\sK}{{\mathcal K}}
\newcommand{\sL}{{\mathcal L}}
\newcommand{\sM}{{\mathcal M}}
\newcommand{\sN}{{\mathcal N}}
\newcommand{\sO}{{\mathcal O}}
\newcommand{\sP}{{\mathcal P}}
\newcommand{\sQ}{{\mathcal Q}}
\newcommand{\sR}{{\mathcal R}}
\newcommand{\sS}{{\mathcal S}}
\newcommand{\sT}{{\mathcal T}}
\newcommand{\sU}{{\mathcal U}}
\newcommand{\sV}{{\mathcal V}}
\newcommand{\sW}{{\mathcal W}}
\newcommand{\sX}{{\mathcal X}}
\newcommand{\sY}{{\mathcal Y}}
\newcommand{\sZ}{{\mathcal Z}}
\newcommand{\A}{{\mathbb A}}
\newcommand{\B}{{\mathbb B}}
\newcommand{\C}{{\mathbb C}}
\newcommand{\D}{{\mathbb D}}
\newcommand{\E}{{\mathbb E}}
\newcommand{\F}{{\mathbb F}}
\newcommand{\G}{{\mathbb G}}
\newcommand{\HH}{{\mathbb H}}
\newcommand{\I}{{\mathbb I}}
\newcommand{\J}{{\mathbb J}}
\newcommand{\M}{{\mathbb M}}
\newcommand{\N}{{\mathbb N}}
\renewcommand{\P}{{\mathbb P}}
\newcommand{\Q}{{\mathbb Q}}
\newcommand{\R}{{\mathbb R}}
\newcommand{\T}{{\mathbb T}}
\newcommand{\U}{{\mathbb U}}
\newcommand{\V}{{\mathbb V}}
\newcommand{\W}{{\mathbb W}}
\newcommand{\X}{{\mathbb X}}
\newcommand{\Y}{{\mathbb Y}}
\newcommand{\Z}{{\mathbb Z}}
\title[Isotrivial families of elliptic surfaces]{
On the isotriviality of families of elliptic surfaces}
\author[Keiji Oguiso]{Keiji Oguiso${}^*$}
\address{Math. Sci., University of Tokyo, 153-8914 Komaba
Meguro Tokyo, Japan}
\email{keiji.oguiso@uni-essen.de /
oguiso@ms.u-tokyo.ac.jp}
\thanks{\hspace{-.26cm}${}^*$Supported by a fellowship of the Humboldt
foundation}\author{Eckart Viehweg}
\address{Universit\"at GH Essen, FB6 Mathematik, 45117 Essen, Germany}
\email{ viehweg@uni-essen.de}
\thanks{This work has been partly supported by the DFG
Forschergruppe ``Arithmetik und Geometrie''}
\maketitle
A family $f: X \to B$ of projective complex manifolds is called
birationally isotrivial, if there exists a finite cover $B' \to
B$, a manifold $F$ and a birational map $\varphi $ from $F
\times B'$ to $X \times_B B'$. The morphism $f$ is isotrivial,
if $\varphi$ can be chosen to be biregular.

One can ask, tempted by the corresponding property
for families of curves, whether $f$ is birationally
isotrivial whenever $B$ is an elliptic curve or $\C^*$ and the
Kodaira dimension of a general fibre non-negative. Assuming that
all fibres of $f$ are minimal models, one could even hope that
$f$ is isotrivial.

Both problems have an affirmative answer, if local Torelli
theorems hold true for the fibres of $f$ (or, as explained in
\ref{variations}, for some \'etale cover), and both have been
solved by Migliorini \cite{Mig} and Kov\'acs \cite{Kov2} for
families of surfaces of general type (see also \cite{Zha},
\cite{dCa} or \cite{BV}).
In this note we want to extend their methods to surfaces of
Kodaira dimension one and thereby complete the proof of the
following theorem.

\begin{thm} \label{main-theorem}
All smooth projective families of minimal surfaces of
non-negative Kodaira dimension over complex elliptic curves or
over $\C^*$ are isotrivial.
\end{thm}

The projectivity assumption is essential. Indeed there exist
smooth, highly non-projective families of $K3$-surfaces over
$\P^1$, called twistor spaces.

Let $M_h$ be the quasi-projective moduli scheme of polarized
manifolds with numerically effective canonical divisor and
Hilbert polynomial $h$ (see \cite{Vie1}). If $Y$ is a complex
algebraic manifold, $\Phi:Y \to M_h$ a morphism, \'etale over its image,
and if $\Phi$ is induced by a ``universal'' family, then
\ref{main-theorem} implies that $Y$ is algebraically hyperbolic
for $\deg (h) =2$ (see also \cite{Kov3}).

If $\bar{Y}$ is a smooth compactification with $S = \bar{Y} - Y$
a normal crossing divisor, one might hope, that
$\Omega^{1}_{\bar{Y}} (\log S)$ (or some symmetric product)
contains a subbundle $\sF$, isomorphic to $\Omega^{1}_{\bar{Y}}
(\log S)$ over $Y$, with $\sF$ numerically effective and $\det( \sF)$ big.
This positivity property holds true for moduli schemes of
curves, and it has recently been verified by Zuo \cite{Zuo}
if the fibres of the universal family over $Y$ satisfy the local
Torelli theorem.

If $B$ is an elliptic curve, or if the fibres $X_b$ of $f$ allow
an \'etale cover which is an elliptic surface without multiple
fibres, the proof of the isotriviality is quite easy. In the
first case, the proof is given at the beginning of section~4,
in the second the necessary arguments are sketched in
\ref{no-multiple} and \ref{easy}, as special cases of the proof of
\ref{main-theorem} for elliptic surfaces, given in section~7.

We thank Egor Bedulev, Fabrizio Catanese, Daniel Huybrechts,
Yujiro Kawamata and Qi Zhang for helpful remarks and comments.
The first named author would like to thank the members of the
``Forschergruppe Arithmetik und Geometry'' at the University of
Essen, in particular H\'el\`ene Esnault, for their hospitality
and help.

\begin{notations} \label{notations}
In discrepancy to the introduction $X$ and $B$ will denote
complex projective manifolds of dimension three and one, and $f:
X \to B$ will be a family of surfaces, i.e. a flat projective
morphism with two dimensional connected fibres $X_b=f^{-1}(b)$.
We fix an open dense subscheme $B_0 \subset B$, such that
$$f_0 =f |_{X_0} : X_0 =f^{-1} (B_0) \>>> B_0$$
is smooth, and we write $S=B - B_0$ and $\Delta = f^* (S)$.

We will call $f$ a family of minimal surfaces, if the
non-singular fibres $X_b$, for $b\in B_0$, are minimal models of
non-negative Kodaira dimension, but we will not require $f$ to
be a relative minimal model in a neighborhood of $f^{-1}(S)$.
 
The dualizing sheaves of $B$, $X$ and of $f$ will be denoted by
$\omega_B$, $\omega_X$ and $\omega_{X/B} = \omega_X \otimes f^*
\omega^{-1}_{B}$.

If $D$ is an effective normal crossing
divisor on $X$, $\Omega^{i}_{X} (\log D)=\Omega^{i}_{X} (\log
D_{\rm red})$ denotes the sheaf of logarithmic differential forms.

Starting from section three, the general fibre $F$ of $f$ is
assumed to be a minimal elliptic surface of Kodaira dimension
$\kappa(F) =1$ and starting with section four, we will assume that
$B$ is an elliptic curve and $S = \emptyset$, or that $(B,S) = (
\P^1, \{ 0, \infty \})$.
\end{notations}
\section{Families of surfaces and isotriviality}

The positivity results for direct images of powers of dualizing
sheaves, due to Fujita, Kawamata and the second named author (see
\cite{Mor}, 7.2 and the references given there) can
be presented in a nice form, if the base is a curve and if the
smooth fibres are minimal.

\begin{definition} \label{with-respect}
Let $X$ be a projective manifold and $U \subset X$ an open dense
subset. An invertible sheaf $\sL$ on $X$ is called
\begin{enumerate}
\item[i)] semi-ample with respect to $U$, if for some $\mu_0$
and all multiples $\mu$ of $\mu_0$ the map
$$
\varphi_{\mu} : H^0 (X , \sL^{\mu}) \otimes_{\C} \sO_X \>>>
\sL^{\mu}
$$
is surjective over $U$.
\item[ii)] ample with respect to $U$, if $\sL$ is semi-ample
with respect to $U$ and if $\varphi_{\mu}$ induces an embedding
$U \to \P (H^0 (X, \sL^{\mu}))$ for $\mu$ sufficiently large.
\end{enumerate}
\end{definition}

\begin{lemma} \label{positivity}
Let $f: X \to B$ be a family of minimal surfaces of non-negative
Kodaira dimension, smooth over $B_0 = B - S$.
\begin{enumerate}
\item[a)] Then $f_* \omega^{\nu}_{X/B}$ is numerically
effective, for all $\nu \geq 1$.
\item[b)] If $f$ is semi-stable, then the following conditions
are equivalent:
\begin{enumerate}
\item[i)] For some $\nu_0 > 0$ and for all multiples $\nu$ of
$\nu_0$ $f_* \omega^{\nu}_{X/B}$ is ample.
\item[ii)] There exists some $\eta > 0$ such that $f_*
\omega^{\eta}_{X/B}$ contains an ample subsheaf.
\item[iii)] $\omega_{X/B}$ is semi-ample with respect of $X_0 =
X - f^{-1} (S)$ and for a general fibre $F$ of $f$ one has
$\kappa (\omega_{X/B} ) = \kappa (F) +1$.
\item[iv)] $f$ is not birationally isotrivial.
\end{enumerate}
\end{enumerate}
\end{lemma}

\begin{corollary} \label{isotrivial-cover}
Let $\tau : Y \to X$ be generically finite. If $f \circ \tau : Y
\to B$ is birationally isotrivial, then the same holds true for
$f: X \to B$.
\end{corollary}

\begin{proof}
We may assume both, $f$ and $f \circ \tau$ to be semi-stable.
The natural inclusion $\omega_{X/B} \to \tau_* \omega_{Y/B}$
induces an inclusion $f_* \omega^{\nu}_{X/B} \to (f\circ \tau)_*
\omega^{\nu}_{Y/B}$, for all $\nu > 0$. Hence if $f$ is not
birationally isotrivial, the condition ii) in \ref{positivity}
b) is satisfied.
\end{proof}

For a smooth projective family $f_0 : X_0 \to B_0$ consider
the polarized variation of Hodge-structures $R^2 f_{0*}
\C_{X_0}$. If $B_0$ is an elliptic curve or $\C^*$, then this
variation of Hodge-structures is necessarily trivialized over
some \'etale cover $B'_0 \to B_0$. In fact, the
induced morphism from the universal cover $\C$ of $B_0$
to the period domain of polarized
Hodge-structures is constant (see for example \cite{Sch}, \S 3).
Combined with \ref{isotrivial-cover} one obtains:

\begin{corollary} \label{variations}
If there exists an \'etale covering $\tau _0 : Y_0 \to X_0$,
such that the fibres of $f_0 \circ \tau_0$ satisfy the local
Torelli theorem, and if $B_0$ is an elliptic curve over $\C^*$,
then $f$ is birationally isotrivial.
\end{corollary}

\begin{remark} \label{kappa-0}
The assumptions of \ref{variations} hold true for all families
of minimal surfaces of Kodaira dimension zero. The same argument
can be used to prove the corresponding statement for families of
curves of genus $g \geq 1$.
\end{remark}

For families of minimal surfaces the birational isotriviality is
equivalent to the isotriviality. As well-known, the
trivialization even exists over an \'etale cover of $B_0$.

\begin{lemma} \label{trivial}
A smooth projective family $f_0 : X_0 \to B_0$ of minimal surfaces (or
curves) of non-negative Kodaira dimension is birationally
isotrivial, if and only if there exists a finite \'etale cover
$B'_0 \to B_0$ and a surface (or curve) $F$ with
$$X_0\times_{B_0} B'_0 \simeq F \times B'_0.$$
\end{lemma}

\begin{proof}
It is easy to find a finite cover $B''_0 \to B_0$ and an
isomorphism
$$
\varphi : X_0 \times_{B_0} B''_0 \> \sim >> F
\times B''_0
$$
of polarized manifolds. In fact, there exists a
coarse moduli space $M_h$ of polarized manifolds, and Koll\'ar
and Seshadri constructed a finite cover of $M_h$ which carries a
universal family (see \cite{Vie1}, p. 298). Of course one may
assume $B''_0 \to B_0$ to be Galois with group $G$. In different
terms, one has a lifting of the Galois action on $B''_0$ to
$F\times B''_0$, giving $X_0$ as a quotient. Let $H$ be the
ramification group of a point $b \in B''_0$. Then $H$ acts
trivially on the fibre $F \times \{ b\}$.

On the other hand, the
automorphism group of a polarized manifold of non-negative
Kodaira dimension is finite, hence the action of $H$ on $F
\times B''_0$ must locally be the pullback under $pr_2$
of the action on $B''_0$. Necessarily the same holds true
globally and
$$
X_0\times_{B_0}(B''_0/H)= (F \times B''_0)/ H = F \times (B''_0 / H).
$$
\end{proof}
\section{A vanishing theorem }

As in \cite{Mig}, \cite{Kov2}, \cite{Zha}, \cite{dCa} or \cite{BV} we will
use vanishing theorems for the cohomology of differential forms
with logarithmic poles. However, we have to allow poles along
some divisor $\Pi$, transversal to the elliptic fibration. In order
to find such a divisor, we will be forced to modify $f_0$ and to
allow some additional singular points in the fibres.

\begin{assumption} \label{van-ass}
Let $X, W$ and $B$ be normal proper algebraic varieties of dimension
three, two and one respectively, and let
$$
\begin{TriCDV}
{X}{ \> g >>}{ W}
{\SE f EE}{\SW W h W}
{B}
\end{TriCDV}
$$
be morphisms with connected fibres. Consider an effective
divisor $\Upsilon$ and a prime divisor $\Pi$ on $X$, and
an invertible sheaf $\sL$ on $W$. Let $B_0 = B -S$ be
open and dense in $B$,
$$
X_0 = f^{-1} (B_0) , \ \ \ \ \ \ \ W_0 = h^{-1} (B_0)
$$
and denote by $f_0$, $g_0$, $\Pi_0$, $\sL_0$ and $h_0$ the
restrictions to $X_0$ and $W_0$, respectively. Assume:
\begin{enumerate}
\item[i)] $\Pi_0$ is a section, i.e. $g|_{\Pi_0} : \Pi_0 \to W_0$ is an
isomorphism.
\item[ii)] $X$ is non-singular and $\Delta = f^* (S)$ as well as $\Delta +
\Pi$ are normal crossing divisors.
\item[iii)] $h_0 : W_0 \to B_0$ is smooth.
\item[iv)] $g_0 : X_0 \to W_0$ is a flat family of curves.
\item[v)] $f_0 : X_0 \to B_0$ is smooth outside of a finite subset $T$ of
$X_0$.
\item[vi)] The sheaf $\sL$ is ample with respect
to $W_0$.
\item[vii)] $h_*\sL^\nu \cong f_*(g^*\sL^\nu \otimes \sO_X(-
\nu\cdot\Upsilon)),$ for all $\nu >0$. In particular $\Upsilon$
is supported in $\Delta$. 
\item[viii)] $\deg \omega_B (S) \geq 0$.
\end{enumerate}
\end{assumption}

\begin{definition} \label{defect}\ \\
\begin{enumerate}
\item[a)] For $\iota : X - T \to X$ define
$\Omega^{i}_{X/B} (\log\Delta)^{\sim} = \iota_* 
\Omega^{i}_{X-T/B} (\log\Delta).$
\item[b)] 
$\displaystyle
\Omega^{i}_{X/B} (\log\Delta)' = \im (\Omega^{i}_{X} (\log \Delta) \>>>
\Omega^{i}_{X/B} (\log\Delta)^{\sim} ).
$
\item[c)] We use the same notation for the sheaves of
differential forms with logarithmic poles along $\Pi$:
\begin{gather*}
\Omega^{i}_{X/B} (\log(\Delta+\Pi))^{\sim} = \iota_* \Omega^{i}_{X-T/B} (\log(\Delta+\Pi))
\mbox{ \ \ \ \ and}\\
\Omega^{i}_{X/B} (\log(\Delta+\Pi))' = \im (\Omega^{i}_{X} (\log( \Delta+\Pi)) \to
\Omega^{i}_{X/B} (\log(\Delta+\Pi))^{\sim} )
\end{gather*}
\end{enumerate}
\end{definition}
Since $\Pi$ does not meet the non-smooth locus $T$ of $f_0$, the
sheaf
$$\Omega^{2}_{X/B} (\log(\Delta+\Pi))'$$
is invertible in a neighborhood of $\Pi$ and
\begin{equation}\label{comp-prime}
\Omega^{2}_{X/B} (\log(\Delta+\Pi))'=\Omega^{2}_{X/B} (\log\Delta)'\otimes \sO_X(\Pi).
\end{equation}
By definition one has the exact sequences
\begin{equation} \label{taut1}
0 \>>> f^* \omega_B (S) \>>> \Omega^{1}_{X} (\log( \Delta+\Pi)) \>>> \Omega^{1}_{X/B}
(\log( \Delta+\Pi))' \>>> 0
\end{equation}
\begin{multline} \label{taut2}
0 \>>> f^* \omega_B (S) \otimes \Omega^{1}_{X/B} (\log(\Delta+\Pi))^{\sim} \>>>
\Omega^{2}_{X} (\log( \Delta+\Pi))\>>>\\
\>>> \Omega^{2}_{X/B} (\log( \Delta+\Pi))' \>>> 0.
\end{multline}
The main result of this section is
\begin{proposition} \label{van}
Assuming \ref{van-ass}
\begin{gather*}
H^0 (X, \Omega^{2}_{X/B} (\log( \Delta + \Pi))' \otimes g^* \sL^{-1}
\otimes \sO_X (\Upsilon - \Pi) \otimes f^* \omega_B (S)^{-2} )= \\
H^0 (X, \Omega^{2}_{X/B} (\log \Delta)' \otimes g^* \sL^{-1} \otimes
\sO_X(\Upsilon)\otimes f^* \omega_B (S)^{-2} ) = 0.
\end{gather*}
\end{proposition}
\begin{remark} \label{van2}
If $f$ is semistable, $\Omega^{2}_{X/B} (\log \Delta)^{\sim}
= \omega_{X/B}$ and $\Omega^{2}_{X/B} (\log\Delta)'$ is
a subsheaf, say $\omega'_{X/B}$, of $\omega_{X/B}$. Then
\ref{van} says that
$$
H^0 (X, \omega'_{X/B}(\Upsilon) \otimes f^* \omega_B (S)^{-2} \otimes g^* \sL^{-1}
) =0.
$$
\end{remark}
\noindent
{\it Proof of \ref{van}}. The statement is compatible with blowing up
$W$ and $X$, as long as the centers are contained in $h^{-1} (S)$ and
$f^{-1} (S)$, respectively. In fact, for $\tau : X' \to X$ and
$\Delta' = \tau^* \Delta$
$$
\Omega^{2}_{X/B} (\log\Delta)'\otimes \sO_X(\Upsilon)
= \tau_* (\Omega^{2}_{X'/B} (\log\Delta')'\otimes
\sO_{X'}(\tau^* \Upsilon)).
$$
Blowing up $W$ (and hence $X$) we may assume that $W$ is non-singular.
For $\mu$ sufficiently large, $\sL^{\mu} (-h^* (S)_{\red})$ is
ample with respect to $W_0$.

Hence, blowing up $W$ and replacing $\mu$ by some multiple, we will find an
effective divisor $\Sigma$ in $W$ such that $\sL^{\mu} (- \Sigma)$
is globally generated and big, and such that $\Sigma_{\red} = h^*
(S)_{\red}$. Moreover, if $\eta : W \to \P (H^0 (W, \sL^{\mu} (-\Sigma))$
denotes the induced morphism, we can also assume that there
exists an effective relatively anti-ample exceptional divisor
$E$. Replacing $\sL^{\mu} (-\Sigma)$ by $\sL^{\mu \cdot \nu} (-
\nu \cdot \Sigma - E)$, we may assume finally that $\sL^{\mu} (-
\Sigma)$ is ample. The assumption \ref{van-ass},
vii), implies that $g^*\Sigma \geq \mu\cdot \Upsilon$.

Since $\Pi_0$ is a section, for some $\rho > 0 $ the map
$$
g^* g_* \sO_X (\rho \cdot \Pi ) \>>> \sO_X (\rho \cdot \Pi)
$$
is surjective over $X_0$. After blowing up $X$, one finds an
effective divisor $\Gamma_1$, supported in $\Delta$, with
$$
g^* g_* \sO_X (\rho \cdot \Pi) \onto{} \sO_X (\rho
\cdot \Pi - \Gamma_1)  .
$$
Let $\Sigma_1 , \ldots , \Sigma_r$ be the irreducible components
of $\Sigma$. For all $\nu$, sufficiently large, and for all $\Sigma' =
\sum^{r}_{i=1} \epsilon_i \Sigma_i \geq 0$, with $\epsilon_i \in
\{0,1\}$,
$$
g^*(\sL^{\mu \cdot \nu})\otimes \sO_X(\rho \cdot \Pi - g^* (\nu \Sigma + \Sigma')
- \Gamma_1 )
$$
is big and generated by its global sections. Choosing $\nu$
larger than $\rho$ and larger than the multiplicities of the
components of $g^*(\Sigma_{\red})$ one finds
$\epsilon_1, \ldots ,\epsilon_r$ such that $N=\nu\cdot\mu$
does not divide the multiplicities of the components of
$$\Gamma = g^* (\nu \cdot \Sigma + \Sigma') + \Gamma_1.$$
By construction $\Gamma_{\red} = \Delta_{\red}$,
$\Gamma \geq N\cdot\Upsilon$,
and $g^*(\sL^N)\otimes \sO_X(\rho \cdot \Pi - \Gamma)$ is
globally generated and big. Let us write
$$
\sL' = g^*(\sL)\otimes\sO_X \Bigl(\Pi - \Bigl[ \frac{\Gamma}{N} \Bigr] \Bigr) =
g^*(\sL) \otimes \sO_X \Bigl(\Pi - \Bigl[ \frac{(N- \rho) \cdot
\Pi + \Gamma}{N} \Bigr] \Bigr).
$$

\begin{claim} \label{van3}
For all $m \geq 0$ and for $i + j < 3$,
$$
H^i (X, \Omega^{j}_{X} (\log( \Delta + \Pi )) \otimes {\sL'}^{-1}
\otimes f^* \omega_B (S)^{-m} ) = 0 .
$$
\end{claim}

Before proving \ref{van3}, let us deduce \ref{van}. Using
\ref{van3} and the long exact cohomology sequence induced by
(\ref{taut2}) $\otimes {\sL'}^{-1} \otimes f^* \omega_B (S)^{-2}$
one obtains an embedding of
\begin{gather*}
H^0:=H^0 (X, \Omega^{2}_{X/B} (\log( \Delta + \Pi ))' \otimes g^* \sL^{-1}
\otimes \sO_X \Bigl( - \Pi + \Bigl[ \frac{\Gamma}{N} \Bigr]
\Bigr) \otimes f^* \omega_B (S)^{-2} )\\
= H^0 (X, \Omega^{2}_{X/B} (\log( \Delta + \Pi ))' \otimes {\sL'}^{-1}
\otimes f^* \omega_B (S)^{-2} )
\end{gather*}
into
$$
H^1:=H^1 (X, \Omega^{1}_{X/B} (\log( \Delta + \Pi))^{\sim} \otimes
{\sL'}^{-1} \otimes f^* \omega_B (S)^{-1} ).
$$
Since $\Omega^{1}_{X/B} (\log( \Delta + \Pi ))' \to \Omega^{1}_{X/B} (\log( \Delta
+ \Pi ))^{\sim}$ is surjective outside of a finite set of points,
$H^1$ is a quotient of
$$
{H'}^1:= H^1 (X, \Omega^{1}_{X/B} (\log( \Delta + \Pi ))' \otimes {\sL'}^{-1}
\otimes f^* \omega_B (S)^{-1} ).
$$
Applying \ref{van3}, for $j =1$, $i =1$, to
(\ref{taut1}), one finds an injective map
$${H'}^1 \>>> H^2 (X, {\sL'}^{-1})$$
and, \ref{van3}, for $j = 0$, $i=2$, implies that both groups are zero.
Hence all the groups, ${H'}^1$, $H^1$ and $H^0$, are zero.
Since $\Bigl[\frac{\Gamma}{N}\Bigr]\geq \Upsilon$ one obtains
\ref{van} from $H^0=0$.
\qed \\[.2cm]
{\it Proof of \ref{van3}.} By the choice of $\sL'$ one has
$$
g^*(\sL^N)\otimes\sO_X (\rho \cdot \Pi - \Gamma) = {\sL'}^{N}
\otimes \sO_X (- (N - \rho ) \cdot \Pi - \Gamma')
$$
for $\Gamma' = \Gamma - N \cdot \Bigl[ \frac{\Gamma}{N}
\Bigr]$. Since $N$ does not divide the multiplicities of the
components of $\Gamma$, one finds $\Gamma'_{\red} =
\Delta_{\red}$. The sheaf $g^*(\sL^N)\otimes\sO_X (\rho
\cdot \Pi - \Gamma)$ contains the inverse image of an ample
invertible sheaf on $W$. All this remains true, if we replace
$\sL'$ by $\sL' \otimes f^* \omega_B (S)^m$, and $\sL$ by $\sL
\otimes h^* \omega_B (S)^m$. So we may assume $m$ to be zero.

If $\delta : X \to \P^M$ is the morphism given by the global
sections of the $\nu$-th power of $g^*(\sL^N)\otimes\sO_X (\rho
\cdot \Pi - \Gamma)$, for $\nu$ sufficiently large,
then $\delta |_{X - \Gamma_{\red}} = \delta |_{X_0}$ can at most
contract components of the fibres of $g_0$. In
particular the maximal fibre dimension of $\delta |_{X_0}$ is
one.

The divisor $H$ of a general section of $g^*(\sL^{N\cdot\nu})\otimes\sO_X
(\rho \cdot \Pi - \Gamma)^{\nu}$ is smooth, and $\Pi +
\Gamma + H$ a normal crossing divisor. By \cite{EV}, 6.2 a) and 4.11 b),
$$
H^i (X, \Omega^{j}_{X} (\log( \Pi + \Gamma + H)) \otimes {\sL'}^{-1}) =0
$$
for $i + j \neq 3$, and
$$
H^i (H, \Omega^{j}_{H} (\log(\Pi + \Gamma) |_{H}) \otimes {\sL'}^{-1})
= 0,
$$
for $i + j \neq 2$. Considering the long exact sequence for
\begin{multline*}
0 \>>> \Omega^{j}_{X} (\log( \Pi + \Gamma)) \otimes {\sL'}^{-1}
 \>>> \Omega^{j}_{X} (\log(\Pi + \Gamma + H)) \otimes
{\sL'}^{-1} \\
\>>> \Omega^{j-1}_{H} (\log( (\Pi + \Gamma)|_{H} ) \otimes
{\sL'}^{-1} \>>> 0
\end{multline*}
one obtains \ref{van3}. \qed

\section{Families of elliptic surfaces}

Let us return to the family $f: X \to B$ of minimal elliptic
surfaces of Kodaira dimension one, with $f_0 : X_0 \to B_0$
smooth. By \cite{Lev} or \cite{Nak}, for all $\nu \geq 0$ and $b
\in B_0$ with $X_b = f^{-1} (b)$,
\begin{equation} \label{base-change}
f_* \omega^{\nu}_{X/B} \otimes \C (b) = H^0 (X_b ,
\omega^{\nu}_{X_b}).
\end{equation}
Fibrewise, for $\nu$ sufficiently large and divisible, $H^0 (X_b
, \omega^{\nu}_{X_b})$ defines the Iitaka map $X_b \to W_b$ to a
non-singular curve $W_b$, and by (\ref{base-change})
$$
f^* f_* \omega^{\nu}_{X/B} \>>> \omega^{\nu}_{X/B}
$$
defines the relative Iitaka map
$$
X \onto{g} W \subset \P (f_* \omega^{\nu}_{X/B} ),
$$
whose restriction $g_0$ to $X_0$ is a morphism, and $W_0 = g
(X_0)$ is smooth over $B_0$.

Blowing up $X$, as always with centers in $f^{-1}(S)$, we can factor $f$ as
$$
\begin{TriCDV}
{X}{\> g >>}{W}
{\SE f EE}{\SW W h W}
{B}
\end{TriCDV}
$$
where $X$, $W$ and $B$ are non-singular projective manifolds and
where $g_0 : X_0 \to W_0$ is a flat projective family of curves.
Using the corresponding property for $X_b \to W_b$, one finds
$$
\dim H^i (g^{-1} (w), \omega^{\nu}_{g^{-1} (w)} ) = 1,
$$
for $i = 0 , 1$ and $w \in W_0$.
Hence $g_{0*} \omega^{\nu}_{X_0/W_0}$ is invertible for all $\nu
\geq 0$. Moreover,
$$g^{*}_{0} g_{0*} \omega_{X_0/W_0} =
\omega_{X_0/W_0} (- \tilde{\Gamma}^{(0)})$$ for some divisor
$\tilde{\Gamma}^{(0)}$.

We will need several properties of elliptic threefolds,
i.e. threefolds with an elliptic fibration. The results needed, due to
Kawamata, Fujita, Nakayama, Miranda, Dolgachev-Gross and Gross are recalled in
\cite{Gro}, together with more precise references. For elliptic
threefolds occurring as the total space of a family of elliptic surfaces,
\cite{FM} is an excellent source. The properties and definitions
needed from the theory of elliptic surfaces, in particular
Kodaira's classification of the singular fibres, can be found in \cite{BPV}.

By \cite{Gro}, Lemma 1.2, blowing up $W$ with centers in $W-
W_0$ one finds a flat relative minimal model $g_m : X_m \to W$,
extending $g_0 : X_0 \to W_0$. Before stating this result in
\ref{relative-mm}, we will use it to define the multiple locus
and the discriminant divisor. In fact to this aim it would be
sufficient to know the existence of $g_m$ over a subscheme $W_1$
with ${\rm codim}(W-W_1) \geq 2$.

Let $\Delta (g_m)$ be the 
smallest subvariety such that
$$g_m^{-1} (W - \Delta (g_m)) \>>> W -
\Delta (g_m)$$ is smooth.

\begin{notations} \label{type}
An irreducible one-dimensional component of $\Delta(g_m)$ belongs
to one of the following, according to the fibre $E=g_m^{-1}(w)$
over the general point $w$ of the component:
\begin{enumerate}
\item[a)] $E$ is a multiple fibre. We denote those components by
$\Sigma_1, \ldots , \Sigma_r$ and call $\Sigma = \sum^{r}_{i=1}
\Sigma_i$ the multiple locus. To $\Sigma_i$ we attach the
multiplicity $m_i$ of the general fibre, and
$\Gamma_i=g_m^{-1}(\Sigma_i)_{\rm red}$, hence
$m_i \cdot \Gamma_i=g_m^{-1}(\Sigma_i)$.
\item[b)] Let $j : W \to \P^1$ denote the rational map, induced
by the $j$-invariant. Let $D_1, \ldots , D_s$ be the components
of the discriminant locus whose image is $\infty$. To $D_i$ we
attach the multiplicity $b_i$ of $D_i$ in $j^{-1} (\infty)$. In
particular, if the general fibre over $D_i$ is a Newton
polygon, then $b_i$ is nothing but the length of the polygon,
(i.e. type $I_{b_i}$). We write $J_{\infty} = \sum^{s}_{i=1} b_i
D_i$.
\item[c)] If $E$ is not a multiple fibre, nor a Newton polygon,
we denote the corresponding components by $D_{s+1}, \ldots ,
D_{\ell}$ and we attach a number $b_i$ to $D_i$ according to
Kodaira's classification (see \cite{BPV}, for example): \\
\ \\
\begin{tabular}{c|c|c|c|c|c|c|c}
type \ \ \ & \ \ \ I$^{*}_{n}$ \ \ \ & \ \ \ II \ \ \ & \ \ \ III \ \ \
& \ \ \ IV \ \ \ & \ \ \ II$^*$ \ \ \ & \ \ \ III$^*$ \ \ \ & \
\ \ IV$^*$ \ \ \ \\ 
\hline
$b_i$ & 6 & 2 & 3 & 4 & 10 & 9 & 8 \\
\end{tabular}\\[.2cm]
\item[d)] $D = \bigcup^{\ell}_{i=1} D_i$ will be called the
discriminant locus, and
$$\sum^{\ell}_{i=1} b_i D_i=J_\infty+\sum^{\ell}_{i=s+1} b_i
D_i$$ the discriminant divisor.
\end{enumerate}
\end{notations}
Remark that $\Sigma$ and $J_\infty$ can have common components,
corresponding to ${}_{m}I_{n}$. The component of the discriminant
locus with general fibre of type $I_{n}^*$ will occur in
$J_\infty=\sum^{s}_{i=1} b_i D_i$ with multiplicity $n$ and in
$\sum^{\ell}_{i=s+1} b_i D_i$ with multiplicity $6$.

\begin{lemma} \label{relative-mm}
Blowing up $W$ with centers in $W - W_0$, there exists a flat
morphism $g_m : X_m \to W$, with $g^{-1}_{m} (W_0) = X_0$
and $g_m |_{X_0} = g_0$, such that
\begin{enumerate}
\item[a)] $W- W_0$ is a normal crossing divisor.
\item[b)] $X_m$ has at most $\Q$-factorial terminal
singularities.
\item[c)] $g_{m*} \omega_{X_m/W}$ is an invertible sheaf
$\delta$.
\item[d)] $\delta^{12} \simeq \sO_W (\sum^{\ell}_{i=1} b_i D_i) =
\sO_W (J_{\infty} + \sum^{\ell}_{i=s + 1} b_i D_i)$, and, for
all $\nu \geq 0$
$$\omega^{[\nu]}_{X_m/W} = g^{*}_{m} \delta^\nu \otimes \sO_{X_m}
\Bigl( \sum^{r}_{i=1} \frac{\nu (m_i-1)}{m_i} \Gamma_i \Bigr).$$
\item[e)] The $j$-invariant defines a rational map $j : W \to
\P^1$, regular in a neighborhood of $h^{-1} (S)$, and $j^*
(\infty) = J_{\infty}$.
\end{enumerate}
\end{lemma}

By \cite{Gro}, lemma 1.2, \ref{relative-mm} holds true if the
discriminant locus is a normal crossing divisor and if one
allows further blow ups. Hence one obtains \ref{relative-mm}
over the complement in $W$ of finitely many points of $W_0$.
Since $X_0$ is non-singular, since $g_0 : X_0 \to W_0$ is flat and since
$g_{0*}\omega_{X_0/W_0}^\nu$ is invertible, 
b), c) and d) extend to $W_0$.

Let us recall the following property of the multiple locus
$\Sigma$, first observed by Iitaka.

\begin{lemma} \label{multiple}
Keeping the notations introduced above, $\Sigma^{(0)} = \Sigma \cap
W_0$ is \'etale over $B_0$ and the fibres of $(g^{-1}_{0}
\Sigma^{(0)})_{\red} \to \Sigma^{(0)}$ are reduced.
\end{lemma}

\begin{proof}
Let us write again $\Gamma^{(0)} = g^{*}_{0} \Sigma^{(0)}$. Then
$\omega_{X_0} = g^{*}_{0} (g_{0*} \omega_{X_0}) \otimes
\sO_{X_0} (\Gamma^{(0)} - \Gamma^{(0)}_{\red})$.

If $\Sigma^{(0)} \to B_0$ is not \'etale, there exists some $b
\in B_0$ such that $\Sigma^{(0)} |_{W_b}$ contains a multiple point.
This remains true, if we replace $B_0$ by any finite cover $B'_0
\to B_0$. In particular, in order to prove the first part of
\ref{multiple}, we may assume that $\Sigma^{(0)} = \sum^{s}_{i=1}
\Sigma_i$, for $\Sigma_i$ the image of a section of $W_0 \to
B_0$. The same can be assumed for the second part. In fact, if
$\Sigma^{(0)} \to B_0$ is \'etale but some fibre of $\Gamma^{(0)}_{\red}
\to \Sigma^{(0)}$ non reduced, then the same remains true after
replacing $B_0$ by an \'etale covering.

Consider for some $r \geq 1$ a point $v\in W_b$ which lies
exactly on $r$ of the components $\Sigma_i$ of $\Sigma^{(0)}$, say
$$
v \in \Sigma_1 \cap \ldots \cap \Sigma_r \cap W_b.
$$
Let $E$ denote the reduced fibre of $g_b$ or $g$ over $v$, let
$\Gamma_i = (g^* \Sigma_i)_{\red}$ and let $m_i$ be the multiplicity of
$\Gamma_i$ in $g^* \Sigma_i$. Finally let $M$ be the
multiplicity of $E$ as a fibre of $g_b : X_b \to W_b$ and
$\Gamma_i . W_b$ the intersection cycle, a positive multiple of $E$.

For all $\mu \geq 1$ the natural map $g^{*}_{0} g_{0*}
\omega^{\mu}_{X_0} \to \omega^{\mu}_{X_0}$ induces an
isomorphism
$$
g^{*}_{0} g_{0*} \omega^{\mu}_{X_0} \> \cong >>
\omega^{\mu}_{X_0} \Bigl( - \sum^{s}_{i=1} m_i \Bigl< \frac{\mu \cdot (m_i
-1)}{m_i} \Bigr> \Gamma_i \Bigr)
$$
where $\bigl< a \bigr> = a - [a]$ denotes the fractional part of a real
number $a$. Since a similar equation holds true for $g_b$, one
obtains
\begin{equation} \label{equation1}
\sum^{r}_{i=1} m_i \Bigl< \frac{\mu \cdot (m_i -1)}{m_i} \Bigr>
\cdot (\Gamma_i . W_b) = M \Bigl< \frac{\mu (M-1)}{M}
\Bigr> \cdot E.
\end{equation}
Choosing for $\mu$ the lowest common multiple $l = {\rm lcm} (m_1,
\ldots , m_r )$ the left hand side of (\ref{equation1}) is zero,
hence $M$ divides $l$. Choosing $\mu = M$,
one finds that each $m_i$ divides $M$, hence $M = l = {\rm lcm} (m_1,
\ldots , m_r)$. For $\mu = M -1 = r_i \cdot m_i -1$ one has \\
\ \\
$\displaystyle \frac{(M-1) (m_i -1)}{m_i} = r_i \cdot m_i - (r_i +1) +
\frac{1}{m_i} $ \
\ \ and \ \ \
$\displaystyle \frac{(M-1)^2    }{M} = M - 2 + \frac{1}{M}.$ \\
\ \\
Therefore (\ref{equation1}) implies that $\sum^{r}_{i=1} \Gamma_i
. W_b = E$. This is only possible for $r = 1$ and if
$\Gamma_1 . W_b$ is reduced.
\end{proof}

\begin{remark}\label{multiple2}
Let $\Sigma_1$ be an irreducible component of the multiple locus
$\Sigma$ and let $\Gamma_1 = g^{-1} (\Sigma_1)_{\red}$. The
fibres of $\Gamma_1 \cap X_0 \to \Sigma_1 \cap W_0$ are either
smooth elliptic curves or Newton polygons. Assume the latter, i.e. that
$\Sigma_1$ is contained in the discriminant locus. Then
$\Gamma_1$ is non-normal. However, since the fibres of
$\Gamma_1$ over points in $\Sigma_1 \cap W_0$ have at most ordinary double
points as singularities, the non-normal locus must be \'etale
over $\Sigma_1 \cap W_0$, hence over $B_0$. Altogether,
replacing $B_0$ by an \'etale covering, we can assume that
$\Sigma^{(0)} = \Sigma \cap W_0$ consists of sections and that the
same holds true for the non-normal locus of the reduced multiple
divisors.
\end{remark}
In order to apply the vanishing stated in \ref{van}, we would
like to restrict ourselves to semistable families $X\to B$.
However, in doing so, one would have to allow $W$ to be singular,
and \ref{relative-mm} would not apply. The
following technical construction will serve as a replacement.
\begin{lemma} \label{prepared}
Let $f: X \to B$ be a family of elliptic surfaces of
Kodaira dimension one, with $f_0 : X_0 \to B_0=B-S$ smooth and
relatively minimal. Assume
that $S$ consists of at least two points, if $B = \P^1$.
Then there exists a finite covering $\tau : B' \to B$, with $B'_0 =
\tau^{-1} (B_0)$ \'etale over $B_0$ and a diagram
of projective morphisms
$$
\begin{CD}
X' \> \eta' >> X^s \> \sigma' >> X\\
\V g' VV \V V g^s V \V V g V\\
W' \> \eta >> W^s \> \sigma >> W\\
\V h' VV \V V h^s V \V V h V\\
B' \> = >> B' \> \tau >> B
\end{CD}
$$
with (as always, the index ${}_0$ refers to the restrictions to
$B_0$ and $B'_0$):
\begin{enumerate}
\item[i)] $\eta_0$ and $\eta'_0$ are isomorphisms. $\sigma_0$ and
$\sigma'_0$ are fibre products.
\item[ii)] $X'$, $W'$ and $X^s$ are non-singular, $W^s$
is normal with at most rational Gorenstein singularities.
\item[iii)] $f^s=h^s\circ g^s: X^s \to B'$ is semistable,
hence the fibres of $h^s:W^s \to B'$ are reduced, and
for $f'=h\circ g$ the fibres $\Delta'={f'}^{-1} (B'-B'_0)$ and
${h'}^{-1} (B'-B'_0)$ are
normal crossing divisors.
\item[iv)] Let $\Sigma'$ be the multiple locus for $g'$ in
$W'$. Then $\Sigma' \cap W_0$ is the disjoint union of sections,
as well as the non-normal locus of ${g'}^{*} (\Sigma')_{\red} \cap X^0$.
\item[v)] $\delta' = g'_* \omega_{X'/W'}$ is invertible, and
$j : W^s \to \P^1$ is regular in a neighborhood of ${(h^s)}^{-1}
(B'-B'_0)$.
\item[vi)] Let $D'$ denote the discriminant locus. Then ${h'}^{-1}
(B'-B'_0) + D' + \Sigma'$ is a normal crossing divisor in a
neighborhood of ${h'}^{-1}(B'-B'_0)$.
\item[vii)] ${\delta'}^{12} = \sO_{W'} (\sum^{\ell}_{i=1} b_i D'_i ) =
\sO_{W'} (J'_{\infty} + \sum^{\ell}_{i=s+1} b_i D'_i)$, where
$\sum^{\ell}_{i=1} b_i D'_i$ is the discriminant divisor, defined
in \ref{type} (in particular, components corresponding to
$I^{*}_{b}$, occur twice).
\item[viii)] Let $\Sigma'_1, \ldots ,\Sigma'_r$ be the components
of the multiple locus which dominate $B'$. Then for all $\nu >0$
one has
$$
f'_* \omega^{\nu}_{X'/B'} = h'_* \Bigl(\omega^{\nu}_{W'/B'} \otimes
\sO_{W'} \Bigl( \sum^{r}_{i=1} \Bigl[ \frac{\nu \cdot (m_i
-1)}{m_i} \Bigr] \Sigma'_i \Bigr) \otimes {\delta'}^{\nu} \Bigr).
$$
\item[ix)] Let $D'_{s +1} , \ldots , D'_{\ell'}$ be those
components of $\sum^{\ell}_{i=s+1} D'_i$, which dominate $B'$. Then
for all multiples $\nu$ of 12
$$
f'_* \omega^{\nu}_{X'/B'} = h'_* \Bigl(\omega^{\nu}_{W'/B'} \otimes
\sO_{W'} \Bigl( \sum^{r}_{i=1} \Bigl[ \frac{\nu \cdot (m_i
-1)}{m_i} \Bigr] \Sigma'_i
+ \frac{\nu}{12} J'_{\infty} + \sum^{\ell'}_{i=s+1} \frac{\nu \cdot
b_i}{12} D'_i \Bigr) \Bigr).
$$
\item[x)] $g^s_*\omega_{X^s/B'}^\nu =
(\eta\circ g')_* \Omega^2_{X'/B'}(\log \Delta')^\nu$ and both
sheaves are reflexive.
\end{enumerate}
\end{lemma}

\begin{proof} We may assume, that $\Delta+D+\Sigma$ is a normal
crossing divisor and that the $j$-invariant defines a
morphism in a neighborhood of $h^{-1}(S)$.

We choose $B'$ to be ramified over $S$ of order divisible by the
multiplicities of the components of $h^{-1}(S)$, and such that
$X\times_BB'$ has a stable reduction $f^s:X^s\to B'$.
\ref{multiple} and \ref{multiple2} allow to assume that iv) holds true.

Choosing for $W^s$ the normalization of $W\times_BB'$, the
fibres of $h^s$ are reduced and $W^s$ has at most rational
Gorenstein singularities. Obviously $f^s$ factors through
$W^s$.

$W'$ is a desingularization of $W^s$, such
that vi) holds true, and such that the flat relative minimal
model, described in \ref{relative-mm}, exists over $W'$.
If we take for $X'$ any desingularization of this minimal model,
$g'_*\omega_{X'/B'}$ is invertible, and vii) holds true.

Up to now, we obtained the first seven properties, and we
remark, that to this aim, we can replace $B'$ by any larger
covering. The last three properties will follow from the first ones.

Let $D$ be an irreducible component of $(h^s)^{-1}(B'-B'_0)$.
Since $f^s$ is semistable, $(g^s)^{-1}(D)$ must be a reduced
normal crossing divisor. In particular, the proper transform
of $D$ in $W'$ can neither belong to the multiple locus, nor to
the discriminant locus, except perhaps to the part corresponding to
Newton polygons. In particular $D$ will not be a component of
$\sum^\ell_{i=s+1} D'_i$.

Let $(g^s_*\omega_{X^s/B'}^\nu)^\vee$ be the reflexive hull.
If $12$ divides $\nu$ then in a neighborhood of a singular point
$w$ of $W^s$ the sheaf $(g^s_*\omega_{X^s/B'}^\nu)^\vee$
is isomorphic to
$$\omega_{W^s/B'}\otimes\sO_{W^s}(\frac{\nu}{12}j^*(\infty)),$$
where $j:W^s \to \P^1$ is the $j$-invariant.
In fact, $w$ can not lie on transversal components of the multiple or
discriminant locus, and as remarked above, all others are part of
$j^{*}(\infty)$. From property vii)
we obtain an injection
\begin{equation}\label{reflexive}
\eta^* ((g^s_*\omega_{X^s/B'}^\nu)^\vee) \> \subset >>
g'_*\omega_{X'/B'}^\nu,
\end{equation}
hence $g^s_*\omega_{X^s/B'}^\nu=(g^s\circ
\eta')_*\omega_{X'/B'}^\nu = \eta_*g'_*\omega_{X'/B'}^\nu$ is
invertible for all multiples $\nu$ of $12$. Moreover,
since the parts of the multiple locus or of the discriminant
locus, which are missing in the formula ix), are all exceptional components
for $\eta$, we obtain ix), as well.

Property viii) follows from ix). For the equality in x) one just
has to remark that for the fibre $\Delta^s$ of $f^s$ over
$B'-B'_0$ one has
$$
\eta'_* \Omega^2_{X'/B'}(\log \Delta') =
\Omega^2_{X^s/B'}(\log \Delta^s) = \omega_{X^s/B'} = \eta'_*
\omega_{X'/B'}.
$$
Let $\sigma$ be a local section of $(g^s_*\omega_{X^s/B'}^\nu)^\vee$
in a neighborhood of a point of $W^s$, which is blown up in $W'$. By
(\ref{reflexive}) the $12$-th power of this section
is the direct image of a section of
$g'_*\omega_{X'/B'}^{12\nu}$, hence of $(\omega_{W'/B'}\otimes
\delta')^{12\nu}\otimes \sO_{W'}(E)$ with $E\geq 0$ exceptional.
Since ${\delta'}^{12}$ contains the inverse image of an invertible
sheaf  on $W^s$, $\sigma$ must be the direct image of a section of
$(\omega_{W'/B'}\otimes \delta')^{\nu}$ and we obtain the
reflexivity in x) for all $\nu$.
\end{proof}
\begin{remark}\label{prepared2}
Given a covering $B''\to B'$, \'etale over $B_0$, we can assume
in \ref{prepared} that $B'$ dominates $B''$. In fact, in the
proof of \ref{prepared} we just used that iv) holds true, and
that the ramification orders are large enough.
\end{remark}

\section{The proof of \ref{main-theorem} in some special cases
and the Jacobian fibration}

Let $f: X \to B$ be a family of minimal elliptic surfaces of
Kodaira dimension one, with $f_0 : X_0 \to B_0$ smooth, and let
$X \> g >> W \> h >> B$ be the factorization constructed in \S
\ 3.\\[.2cm]
{\it Proof of \ref{main-theorem} for smooth families of elliptic
surfaces of general type over elliptic curves.}
If $B=B_0$ is an elliptic curve, the total space $X=X_0$ of a
family of minimal elliptic surfaces is itself a minimal model,
and the proof of the isotriviality is similar to the one
given in \cite{Mig} for families of surfaces of general type.
 
As in \ref{variations} the polarized variations of Hodge
structures $R^i f_* \C_{X}$ are trivial, hence $R^i f_* \sO_X$
is a free sheaf of degree zero, and by the Leray spectral
sequence and by the Riemann Roch theorem on $B$ and on $X$ one
obtains
\begin{multline*}
- \frac{c_1 (X) . c_2 (X)}{12} = \chi (\sO_X) =
\sum^{2}_{i=1} (-1)^i \chi (R^i f_* \sO_X) =\\
\sum^{2}_{i=1} (-1)^i \deg(R^i f_* \sO_X) = 0.
\end{multline*}
Assume that $f$ is non-isotrivial and let $g : X \to W$ be the
relative Iitaka map. \ref{positivity} implies that
$\omega_{X/B}$ is numerically effective of Kodaira-dimension 2,
and by the canonical bundle formula, for $\nu$ sufficiently
large and divisible, $\omega^{\nu}_{X/Y} = g^* \sA$, with $\sA$
ample on $W$, and $(g^* c_1 (\sA)) . c_2 (X) =0$. For a fibre $W_b$
of $h$, one finds
\begin{equation}\label{inters}
(g^* c_1 (\sA(-W_b))) . c_2 (X) +
(g^* W_b) . c_2 (X)=0.
\end{equation}
On the other hand, since $X$ is a minimal model, \cite{Miy},
3.2, implies that $c_2(X)$ is pseudo-effective. So, choosing
$\nu$ large enough, none of the summands in (\ref{inters}) can be negative.
Thus $c_2(X_b) = (g^* W_b). c_2(X) =0$,  showing that the only
singular fibres of $X_b \to W_b$ are multiple fibres. One obtains
$$
K_{X/B} = g^* \Bigl( K_{W/B} + \sum^{r}_{i=1} \frac{m_i -1}{m_i}
\Sigma_i \Bigr) ,
$$
as $\Q$-divisors. 

By \ref{trivial}, applied to $h : W \to B$, we may assume that
$W = C \times B$ and $h = pr_2$, if $g (W_b) \geq 1$. The same
holds true for $W_b = \P^1$, since the $g_b : X_b \to W_b$ has
at least three multiple fibres in that case. If $g(C) \neq 1$,
for all $i$ the images $pr_1 (\Sigma_i)$ are points,
contradicting the ampleness of the $\Q$-divisor
$$K_{W/B} + \sum^{r}_{i=1}
\frac{m_i -1}{m_i} \Sigma_i.
$$
If $g(C) =1$, then $K_W = K_{W/B}
=0$ and $0 = \deg K_{\Sigma_i} = (K_W + \Sigma_i) . \Sigma_i =
(\Sigma_i)^2.$ Hence $(K_{W/B} + \Sigma)^2 =0$, again contradicting the
ampleness. \qed\\[-.1cm]

A relatively minimal elliptic fibration
$\tilde{\gamma}:\tilde{J} \to W$ is called the
Jacobian-fibration of $g$, if the generic fibre of
$\tilde{\gamma}$ is the Jacobian of the generic fibre of $g$. As
explained in \cite{Gro}, 1.4~-~1.6, even if $g:X\to W$ has a
flat relative minimal model (see \ref{relative-mm}), one can not assume
$\tilde{\gamma}$ to be flat. In fact, one has to exclude
the points, where the discriminant locus has non-normal
crossings, and certain types of collision points. Nevertheless,
by \cite{Gro}, 1.6, the canonical bundle formula
$\omega_{\tilde{J}}=\tilde{\gamma}^* (\omega_W\otimes \delta)$
remains true.

Assume that for $g_0 : X_0 \to W_0$ the multiple locus is empty.
Since the same holds true for the fibres $X_b \to W_b =
h^{-1} (b)$, each fibre of $g_0$ has a reduced component, and
$g_0 : X_0 \to W_0$ has local sections over \'etale
neighborhoods of all points. In this
case, we may choose $\tilde{\gamma}_0 : \tilde{J}_0 \to W_0$, to be locally
in the \'etale topology isomorphic to $g_0 : X_0 \to W_0$.
In particular, $\tilde{J}_0$ is non-singular.

The same remains true, if $X_0$ is non-singular, but if
finitely many of the fibres $X_b \to
W_b$ have isolated singularities. 

We choose a desingularization $\sigma: J \to \tilde{J}$ with
$\sigma^{-1}(\tilde{J}_0) \cong \tilde{J}_0$. The induced
family 
$$
\gamma=\tilde{\gamma}\circ\sigma : J\>>> W
$$
will be called a Jacobian fibration of $g$. 

\begin{lemma} \label{compare}
Assume that $g_0 : X_0 \to W_0$ has no multiple fibres and that
$$
X \> g >> W \> h >> B
$$
satisfies the conditions stated in \ref{prepared}, vii) - ix)
(with $B' = B$). Let $J \> \gamma >> W$ be a
Jacobian fibration. Then, using the notations from
\ref{prepared}
$$
\gamma_* \omega^{\nu}_{J/W} = \delta^{\nu} \ \ \mbox{and} \ \
f_* \omega^{\nu}_{X/B} = ( h \circ \gamma)_* \omega^{\nu}_{J/B}
$$
for all $\nu$ divisible by 12.
\end{lemma}

\begin{proof} By the canonical bundle formula \cite{Gro}, 1.6,
$\gamma^*\delta$ is a subsheaf of $\omega_{J/W}$. Hence
$\delta^\nu$ is an invertible subsheaf of $\gamma_*
\omega^{\nu}_{J/W}$, and since both coincide outside of a finite
number of points, they are the same. The second equality follows
from \ref{prepared} viii). 
\end{proof}

\begin{corollary} \label{no-multiple}
\ref{main-theorem} holds true for families of elliptic surfaces
of Kodaira dimension one and without multiple fibres.
\end{corollary}

\begin{proof} For a Jacobian fibration $\gamma:J\to W$
write $\psi=h\circ \gamma:J \to B$. 
By \ref{prepared} and \ref{prepared2} we can find a
covering $B'\to B$, \'etale over $B_0$, such that the conditions
in \ref{prepared} are satisfied for suitable models of both,
$X\times_BB'$ and $J\times_BB'$. We will drop the ${}'$ and
assume $B=B'$.

The family $f:X\to B$ is birational to the semistable family
$f^s:X^s \to B$. If $f$ is not birationally isotrivial,
\ref{positivity} implies that $f_*\omega^{\nu}_{X/B}$ is ample,
and $\omega_{X/B}$ will be semi-ample with respect to $X_0$.
The property viii) in \ref{prepared} implies that,
$\omega_{W/B} \otimes \delta$ is ample with respect to $W_0$.
By \ref{compare}, the same holds true for $\sL=\gamma_*
\omega_{J/B}$. Choose the effective divisor $\Upsilon$, such that
\begin{equation}\label{xi}
\Omega^2_{J/B}(\log \psi^{-1}(S))\cap \gamma^*\sL =
\gamma^*\sL\otimes \sO_J(-\Upsilon).
\end{equation}
The last condition in \ref{prepared} implies that,
for all $\nu >0$,
$$
\psi_*\Omega^2_{J/B}(\log \psi^{-1}(S))^\nu =
\psi_*\omega_{J/B}^\nu,
$$
hence $\Upsilon$ satisfies the condition vii) in \ref{van-ass}.
$J_0 \to B_0$ is smooth, and choosing $\Pi$ as the closure of
the zero-section of $J_0 \to W_0$ the assumptions in
\ref{van-ass} hold true (with $T=\emptyset$). By \ref{van}
$$
H^0 (J, \Omega^{2}_{J/B} (\log \psi^{-1}(S)) \otimes
\gamma^* \sL^{-1} \otimes \sO_J(\Upsilon)) = 0,
$$
contradicting the choice of $\Upsilon$ in (\ref{xi}).
\end{proof}

Using \ref{variations} and some special considerations for the
case that the $j$-invariant is constant along the fibres $W_b$,
one can replace the reference to \ref{van} in the proof of
\ref{no-multiple} by Saito's local Torelli theorem \cite{Sai}.

\begin{corollary} \label{easy} For $B_0=\C^*$,
\ref{main-theorem} holds true if the general fibre $X_b$ is an
elliptic surface of general type, and if the Iitaka map $g_b :
X_b \to W_b$ satisfies one of the following:
\begin{enumerate}
\item[a)] $g(W_b) \geq 1$.
\item[b)] $W_b \cong \P^1$ and $g_b$ has three or more multiple
fibres.
\item[c)] $W_b \cong \P^1$ and $g_b$ has two multiple fibres of
the same multiplicity $m$.
\item[d)] $W_b \cong \P^1$ and $g_b$ has two smooth multiple
fibres of multiplicity larger than $6$.
\end{enumerate}
\end{corollary}
\noindent {\it Sketch of the proof.} We may assume that the transversal
components $\Sigma_1, \ldots ,\Sigma_r$
of the multiple locus are the images of sections of $X_0\to W_0$. 

In a) $W_0 \to B_0$ is an isotrivial family of curves and by
\ref{trivial} we may assume that $W_0 = C \times B_0$.
Then the multiple locus is of the form $\sum^{r}_{i=1} c_i
\times B_0$. If $W_b=\P^1$ we may choose
an isomorphism $W_0\cong \P^1\times B_0$ with $\Sigma_i = c_i
\times B_0$. 

In the first three cases there exist coverings of $C$ or $\P^1$
with exact ramification order $m_i$ over $c_i$ (see \cite{FK},
IV.9.12, for example). 

In case a) or c) it is easy to describe such a covering explicitly:
Replacing $C$ in a) by an \'etale cover of degree two, we may
assume that the multiplicities of the fibres over $c_{2i}\times B_0$ and
over $c_{2i+1}\times B_0$ are $m_{2i}$. By \cite{EV}, 3.15, the
covering obtained by taking the $m_{2i}$-th root out of the
divisor $c_{2i} + (m_{2i}-1)\cdot c_{2i+1}$ is totally ramified
of order $m_{2i}$ over $c_{2i}+c_{2i+1}$, and nowhere else. The
normalization of the fibred product of the coverings obtained,
is the one asked for. In c) one just takes the $m$-th root out
of the divisor $c_1 + (m-1)\cdot c_2$.

Hence in a), b) or c) there exists a covering $W'_0$, 
ramified over $\Sigma_i$ of order $m_i$ and \'etale
over $W_0-\bigcup_{i=1}^r\Sigma_i$. The normalization $X'_0$ of
$X_0\times_{W_0}W'_0$ is \'etale over $X_0$, hence it remains
smooth over $B_0$. The projection to $W'_0$ has no multiple fibres, 
and \ref{easy} follows from \ref{no-multiple} and
\ref{isotrivial-cover}. 

For d) one shows, as indicated in \ref{moduli-curves},
that after replacing $B_0$ by an \'etale cover, a multiple
component $\Sigma_i$ with multiplicity $m_i$ gives rise to a
morphism from $\Sigma_i$ to the moduli scheme of elliptic curves
with level $m_i$-structure. Since the genus of this moduli
scheme is larger than one, for $m_i > 6$, this map
must be constant. Hence $g^{-1} (\Sigma_i)_{\red} \to \Sigma_i$ is
smooth over $\Sigma_i \cap W_0$. Choose a covering $W'_0 \to W_0$,
ramified of order $m_1 \cdot m_2$ along $\Sigma_1 + \Sigma_2$, and nowhere
else. Then the normalization of $X_0 \times_{W_0} W'_0$ is again smooth
over $B_0$, but without multiple fibres. \qed\\[-.1cm]

Although we will reprove \ref{easy} in section~7, using slightly
different coverings $W'_0 \to W_0$, let us
concentrate for a moment on those families, not covered by
\ref{easy}, a), b) or c), i.e. those
with $B_0 = \C^*$, with $W_b = \P^1$ and with one of the
following:\\[.2cm]
\underline{Case I:} There are two multiple fibres of
multiplicities $m_1 \neq m_2$ in $X_b \to W_b = \P^1$.\\[.2cm]
\underline{Case II:} There is one multiple fibre of
multiplicity $m$ in $X_b \to W_b = \P^1$.\\[-.1cm]

In the first case, we will replace $X \to W$ by a
desingularization $X'$ of the pullback $X \times_{X} W' \to W',$
where $W' \to W$ is totally ramified over $\Sigma_1 + \Sigma_2$
of order $M$, divisible by
$m_1$ and $m_2$. Doing so, the morphism $X'_0 \to B_0$ will no
longer be smooth in a finite subset of $X'_0$. A careful
analysis of the geometry of the multiple fibres in section~6
will allow to choose $M$ in such a way, that $X'$ locally
factors through a finite morphism $X' \to X''$, with $X''$
smooth over $B$. This observation will allow to apply \ref{van}
to $X'_0$, along the same lines used to prove \ref{no-multiple}.
The sheaf $\sL$ will correspond to the inverse image of the $\Q$-divisor
$K_{W/B} + \sum^{2}_{i=1} \frac{m_i -1}{m_i} \Sigma_i + \delta$
on $W'$. 

The same construction (with $m = m_1$ and $m_2 =1$) works in case II,
if one is able to choose the second section $\Sigma_2$ in such a
way that it only meets components of the discriminant locus corresponding to
reduced singular fibres (types $I_n$, $II$, $III$ or $IV$).
To find such a section, we will have to study the discriminant
locus in section~5. There we will use in an essential way that
$\chi (\sO_{X_b}) \geq 2$, a condition which fortunately holds
true in case II. 

\section{Constantness of the Weyl system}

If $\chi (\sO_{F}) \geq 2$, for a general fibre $F$ of $f: X \to
B$, then the triviality of the variations of Hodge structures
forces the part of the discriminant locus which corresponds to
non-reduced non-multiple fibres to be \'etale over $B_0$.

\begin{proposition} \label{weyl}
Let $f_0 : X_0 \> g_0 >> W_0 \> h_0 >> B_0$ be a smooth
projective family of minimal elliptic surfaces with $\chi
(\sO_{X_b}) \geq 2$ and $\kappa(X_b) =1$, for all $b \in B_0$ and $X_b =
f^{-1} (b)$. Assume that $B_0 = \C^*$ or that $B_0$ is an
elliptic curve. Let 
$$
D^{(0)} = \sum^{\ell}_{i=s+1} D_i
$$
be the part of the discriminant locus in $W_0$, which corresponds to
singular fibres of types $I^{*}_{j} \ (j \geq 0)$, $II^*$,
$III^*$ or $IV^*$. Then $D^{(0)}$ is \'etale over $B_0$,
the restriction $g_0^{-1}(D^{(0)}) \to D^{(0)}$ is locally
equi-singular, and $D^{(0)} \cap \Sigma^{(0)} = \emptyset$ for
the multiple locus $\Sigma^{(0)}$ of $X_0 \to W_0$.
\end{proposition}

The condition $\chi (\sO_{X_b}) \geq 2$ is needed in the proof
of the following description of the -2 classes in the
N\'eron-Severi group $NS(X_b)$ of $X_b$.

\begin{lemma} \label{neron-severi}
Let $g_b : X_b \to W_b$ be a minimal elliptic surface of Kodaira
dimension one with $\chi (\sO_{X_b}) \geq 2$. Define the numbers
$n_i$, $m_j$ and $l_k$ as the number of reducible fibres,
according to the following list: \\
\ \\ {\rm
\begin{tabular}{l|c|c|c|c|c|c|c}
type & ${}_mI_i$ & \ $IV$ or \ &
\ $III$ or \ & \ \ $I^{*}_{j}$ \ \ & \ \ \ $II^*$ \ \ \ & \ \ $III^*$ \ \
& \ \ $IV^*$ \ \ \\
& $(m \geq 1, i \geq 4)$ & ${}_mI_3$ &
${}_mI_2$ & & & & \\
\hline
number &&&&&&&\\
of fibres & $n_i$ & $n_3$ & $n_2$ & $m_j$ & $l_8$ & $l_7$
& $l_6$\\
\hline
Euler &&&&&&&\\
number & $i$ & $4$ or $3$ & $3$ or $2$ & $j + 6$ & $10$ & $9$ & $8$
\end{tabular}}\\[.3cm]
Let \hfill $N_b = \bigl< \alpha \in NS (X_b); \ (\alpha . F) =0
\ \mbox{and} \ (\alpha . \alpha) = -2 \bigr> / \Q \cdot F \cap
NS(X_b)$ \hspace*{\fill}\\[.1cm] 
be the root lattice. Then
\begin{enumerate}
\item[i)] $N_b$ is generated by the classes of irreducible
components of reducible fibres of $g_b$.
\item[ii)] The numbers $n_i , m_j$ and $l_k$ are uniquely
determined by $N_b$ and by its decomposition
$$
N_b \simeq \bigoplus_{i \geq 2} A^{\oplus n_i}_{i} \oplus
\bigoplus_{j \geq 0} D^{\oplus m_j}_{j+4} \oplus
\bigoplus^{8}_{k=6} E^{\oplus l_k}_{k}
$$
in indecomposable sublattices.
\end{enumerate}
\end{lemma}

\begin{proof}
The assertion ii) follows from i) and from the well-known
uniqueness of the decomposition of $N_b$ (see for example \cite{Hum},
Proposition 11.3).

For i) let $[D] \in NS (X_b)$ be a representative of a class
$\alpha \in N$, with $(D.D) = - 2$. Since $K_{X_b}$ is
numerically equivalent to $a \cdot F$, for some $a \in \Q$, one
obtains from the Riemann Roch formula
$$
\chi (\sO_{X_b} (D))= \frac{D.(D-K_{X_b})}{2} + \chi (\sO_{X_b})
= - 1 + \chi (\sO_{X_b}) > 0.
$$
Therefore $H^0 (\sO_{X_b} (D)) \neq 0$ or $H^0 (\sO_{X_b} (K_{W_b} -
D)) \neq 0$. Since $[K_{W_b} - D] = [-D]$ in $N_b$, replacing
$\alpha$ by $- \alpha$ we may assume that $\alpha$ is
represented by an effective divisor $\Sigma \alpha_i D_i$. 
Since $0 = (\alpha . F) = \Sigma \alpha_i (D_i . F)$ one finds the
$D_i$ to be irreducible components of the fibres of $g_b$, and
one may assume that those are components of reducible fibres.
One obtains i) from the classification of the singular fibres
(see \cite{BPV}, for example).
\end{proof}

\noindent{\it Proof of \ref{weyl}.}
Let $\sB = \C \to B_0$ be the universal covering of $B_0$ and
denote the pullback of $X_0 \to W_0 \to B_0$ by
$$
\tilde{f} : \sX \> \tilde{g} >> \sW \> \tilde{h} >> \sB .
$$
Then $R^2 \tilde{f}_* \Z_{\sX}$ is a constant system, i.e. we
have a global marking
$$
\tau : R^2 \tilde{f}_* \Z_{\sX} \> \cong >> H^2 \times \sB
$$
for $H^2$ a lattice isomorphic to $H^2 (X_b, \Z)$. Recall that
any invertible sheaf $\sL$ on $\sX$ defines a constant subsystem
$$
c_1 ( \sL |_{X_b})_{b \in\sB} \subset R^2 \tilde{f}_* \Z_{\sX}.
$$
In particular, if $\sH$ is the inverse image of a relative ample
invertible sheaf on $X_0 \to B_0$, we can define the constant
system $(R^2 \tilde{f}_*\Z_{\sX} )_{{\rm prim}} = [\sH]^{\perp}$
in $R^2 \tilde{f}_* \Z_{\sX}$. Restricting $\tau $ one obtains
an isomorphism 
$$
\tau^{\perp} : (R^2 \tilde{f}_* \Z_{\sX})_{{\rm prim}} \> \cong
>> H^{\perp} \times \sB ,
$$
where $H^{\perp} \subset H^2$ is a sublattice. Using
$\tau^{\perp}$, we define the global period map
$$
p : \sB \>>> {\rm Grass} (k, H^{\perp}) \mbox{ \ \ \ by \ \ \ } p (b) =
(\tau^{\perp} (H^0 (X_b, \Omega^{2}_{X_b} )) \subset H^{\perp}
\otimes \C ).
$$
Since $\sB = \C$ is a Zariski-open subset of $\P^1$, we may
apply \cite{Sch}, Theorem 7.22, and we find $p$ to be constant.
Hence $\tilde{f}_* \Omega^{2}_{\sX/\sB}$ is a flat
vector bundle, that is, there exists a linear subspace $T
\subseteq H^{\perp} \otimes \C$, such that $\tau^{\perp}$
induces an isomorphism
$$
\tilde{f}_* \Omega^{2}_{\sX/\sB} \> \cong >> T \otimes
\sO_{\sB}.
$$
We define $NS = T^{\perp} \cap H^{2}_{\Z}$ and consider the
corresponding constant system
$$
\tau : \sN \sS \> \cong >> N S  \times \sB.
$$
By the Leftschetz (1,1) Theorem , $\sN \sS$ is the system
consisting fibrewise of the N\'eron-Severi groups $NS(X_b)$ (Note
that in general, the system $NS(X_b)$ is far from being
constant).

Next we consider the constant subsystem $\sC_1$, defined by the
relative dualizing sheaf $\omega_{\sX/\sB}$, and the sublattice
$c_1 \subset H^2$, with $\tau (\sC_1)= c_1 \times \sB$. Of
course, $\sC_1 \subset \sN \sS$ and $(c_1 . c_1) =0$. Taking the
quotient we obtain
$$
\tau' : \sS = \sC^{\perp}_{1} / \sC_1 \otimes_{\Z} \Q \cap \sN
\sS \> \cong >> S \times \sB = c^{\perp}_{1}/c_1 \otimes_{\Z} \Q
\cap NS.
$$
Up to now, we obtained a constant system $\sS$ consisting fibrewise of
$$
\{ a_b \in NS (X_b) ; \ a_b . c_1 (\omega_{X_b}) =0 \}
$$
modulo rational multiples of $c_1 (\omega_{X_b})$. Since $\kappa
(X_b) =1$, $c_1 (\omega_{X_b})$ is some positive multiple of the
class of a fibre $[F_b]$ of $X_b \to W_b$. Note that the
intersection form on $H^2 (X_b , \Z)$ descends to the one on
$\sS$, since the condition
$$
(\alpha . c_1 (\omega_{X_b})) = (\beta . c_1 (\omega_{X_b})) =0
$$
implies that
$$
(\alpha . \beta ) = ((\alpha + a \cdot c_1 (\omega_{X_b})).( \beta
+ b \cdot c_1 (\omega_{X_b}))).
$$
$S$ is a negative definite lattice. Consider the sublattice
$$ N:
= \bigl< \alpha \in S ; (\alpha . \alpha ) = -2\bigr>
$$
and the corresponding constant system
$$
\sN \> \cong >> N \times \sB.
$$
This $\sN$ is a family of lattices, which fibrewise corresponds
to the lattice $N_b$ described in \ref{neron-severi}. In
particular, the numbers $n_i$, $m_j$, and $l_k$ defined in
\ref{neron-severi} are independent of $b \in \sB$, and by
definition of $n_i \ (i \geq 4)$, $m_j$, and $l_k$ one obtains
\begin{claim} \label{fibres}
The number of singular fibres of type $I_i \ (i \geq 4)$,
$I^{*}_{j}$, $II^*$, $III^*$ and $IV^*$ in $X_b \to W_b$ is
independent of $b \in B_0$.
\end{claim}

To finish the proof of \ref{weyl} we need the constantness of
the local Euler numbers. For $p \in \sW$ choose small disks
$\Delta^{2}_{(x,y)} \subset \sW$ with center $p$, and $\Delta_x
\subset \sB$ with center $\tilde{h} (p)$, such that $\tilde{h}
|_{\Delta^{2}_{(x,y)}} : \Delta^{2}_{(x,y)} \to \Delta_x$ is the
projection of the first factor. For $\epsilon > 0$ sufficiently
small and $\alpha \in \C$, with $|\alpha | < \epsilon$, we write
$L_{\alpha, \epsilon} = (x- \alpha \cdot y =t)$, and $\Delta^2 =
\bigcup_{|t|< \delta} L_{\alpha, t}$. Hence $\tilde{g}^{-1}
(L_{\alpha , t}) \to L_{\alpha , t}$ is a family of smooth,
local elliptic surfaces, parameterized by $t$.

The Euler numbers of $\tilde{g}^{-1} (L_{\alpha ,t})$ are
independent of $t$, and they are the sum of the Euler numbers of
the singular fibres of $\tilde{g}^{-1} (L_{\alpha, t}) \to
L_{\alpha , t}$.

Let again $D_1, \ldots , D_\ell$ be the components of the
discriminant locus and let $e_i$ be the Euler number of the
general fibre over $D_i$. We assume, renumbering $D_1,\ldots
,D_\ell$ if necessary, that $e_1 \leq e_2 \leq \ldots \leq e_\ell$.
Assume that for $i_0$ the divisor $\sum^{\ell}_{i=i_0+1} D_i$ is
\'etale over $B_0$, but $\sum^{\ell}_{i=i_0} D_i$ is not. Hence
there exists $p \in \sW$, where for a suitable choice of
$\epsilon$,
$$
\tilde{g}^{-1} (L_{0,0}) \>>> L_{0,0}
$$
has just one singular fibre, whereas the number of singular
fibres in
$$
\tilde{g}^{-1} (L_{0,t}) \>>> L_{0,t}
$$
is larger than or equal to two. Hence the Euler number of
$\tilde{g}^{-1} (p)$ is strictly larger than $e_{i_0}$. By
\ref{fibres} this is only possible for $e_{i_0} < 4$.

Since $e_i \geq 6$, for the components corresponding to singular
fibres of types $I^{*}_{b}$, $II^*$, $III^*$, $IV^*$, we obtain
that the union of the corresponding components is \'etale over
the base. Also, those components can not meet the multiple
locus, since the reduced fibre of an intersection point must be
a Newton polygon of length larger than or equal to $6$,
contradicting again \ref{fibres}. \qed

\begin{remark} \label{non-etale}
The method used to prove \ref{weyl} gives a bit more. The
constantness of the local Euler numbers and \ref{fibres} exclude
for example, that some $p \in W_0$ lies on two components $D_1$
and $D_2$ of the discriminant locus, which correspond to fibres
of type $I_{b_1}, I_{b_2}$ with $b_1 + b_2 \geq 5$.

Nevertheless, the method is not strong enough, to imply the
\'etaleness of the whole discriminant locus. For example it can
not exclude that $D_1$, a component corresponding to $I_1$, has
a cusp. Such examples exist locally, and two $I_1$-fibres
degenerate towards a $II$-fibre in such a point.
\end{remark}

\section{Standard modifications of multiple fibres }

Let $f_0: X_0 \to W_0 \to B_0$ be a smooth family of minimal
elliptic surfaces. We assume that the multiple locus $\Sigma^{(0)} =
\sum^{r}_{i=1} \Sigma_i$ and the non-normal locus of $g^{-1}_{0}
(\Sigma^{(0)})_{\red}$ consists of the union of disjoint sections.
Let $m_1, \ldots , m_r$ be the multiplicities of
$g^{-1}_{0}(\Sigma_1), \ldots , g^{-1}_{0}(\Sigma_r)$,
respectively. The multiple locus can meet other components of
the discriminant locus. An example, due to Moishezon, is given in
\cite{FM}, 7.4, where $\Sigma_i$ is not contained in the
discriminant locus, but meets a component $D_1$ of type $I_n$.
In this example, $\Sigma_i$ is an $m_i$-fold tangent to $D_1$.
In fact, it is easy to show, that $\Sigma_i$ can only meet the
discriminant locus in components of $J_{\infty}$, and this
intersection can not be transversal. For components of
type $I^{*}_{n}$, $II^*$, $III^*$ and $IV^*$, this has been part
of \ref{weyl}, at least if $B_0 = \C^*$ or an elliptic curve.

Consider a finite covering $W'_0 \to W_0$ which is totally
ramified of order $m_i$ and \'etale over $W_0 - \Sigma_i$, in a
neighborhood of $\Sigma_i$. The normalization $X'_0$ of $X_0
\times_{W_0} W'_0$ is \'etale over $X_0$, and one obtains an
\'etale Galois cover $\Gamma'_{i} \to \Gamma_i = g^{-1}_{0}
(\Sigma_i)_{\red}$. Hence $\Gamma'_i$ has a fixed point free
action of $\Z/m_i \Z$, and the only singular fibres of
$\Gamma'_i \to \Sigma'_i = (\Sigma_i \times_{W_0} W'_0)_{\red}$
are smooth elliptic curves or Newton polygons of length divisible by
$m_i$. 
\begin{remark} \label{moduli-curves}
The $j$-invariant might be non constant along $\Sigma_i$.
Although not needed in the sequel, let us point out some obvious
obstructions for this to happen, in case $B_0 = \C^*$. If $J_i$
denotes the Jacobian of $\Gamma'_i$, we can assume (replacing
$B_0$ by an \'etale cover) that $J_i$ has a level
$m_i$-structure. Hence the $j$-invariant factors through
$\Sigma'_i \to X_1 (m_i)$, where $X_1 (m_i)$, is the moduli
curve parameterizing elliptic curves with a level
$m_i$-structure. By \cite{Shi}, 1.6.4, $g(X_1 (m_i)) =0$ implies
that $m_i \leq 6$. Using the additional information, that the
translation by one of the sections of order $m_i$ can only have
fixed points in two fibres (the ones over $\{ 0 , \infty \}
\subseteq \P^1 = B)$, one can exclude the case $m_i = 6$, but
not the others.
\end{remark}
For the proof of \ref{main-theorem} we will need a slightly
different description of the multiple locus. For example, if
$X_b \to W_b$ is an elliptic surface with two multiple fibres of
multiplicities $m_1$ and $m_2$, we have to replace $W_0$ by
a covering $W'_0$, totally ramified of order divisible by ${\rm
lcm}(m_1,m_2)$. 
The normalization of the pullback of $X_0$ will have no multiple
fibres anymore, but it might not allow a model, smooth over $B_0$.
So in order to apply \ref{van}, along the same line we did in
\ref{no-multiple}, we have to construct a model $X'_0$ for
which we control the sheaf $\Omega_{X'_0/B_0}^2(\log
\Delta')'$, defined in \ref{defect}. 

\begin{lemma} \label{multiple-locus}
Assume that the multiple locus $\Sigma^{(0)}$ of
$$
X_0 \> g_0 >> W_0 \> h_0 >> B_0
$$
consists of sections, as well as the non-normal locus of
$g^{-1}_{0} (\Sigma^{(0)})_{\red}$. Then we can attach to each
component $\Sigma_i$ of $\Sigma^{(0)}$ a number $\mu_i$, divisible
by the multiplicity $m_i$ of $\Gamma_i = g^{-1}_{0} (\Sigma_i)$
with the following property.

Let $\tau : W'_0 \to W_0$ be a covering, totally ramified of
order $M$, divisible by $\mu_i$, along $\Sigma_i$ and unramified
over $U - \Sigma_i$ for a neighborhood $U$ of $\Sigma_i$ in
$W_0$. Then there exists a commutative diagram of projective
morphisms
$$
\begin{CD}
X'_0 \> \tau' >> X_0 \\
\V g'_0 VV \V V g_0 V \\
W'_0 \> \tau >> W_0 \\
\V h'_0 VV \V V h_0 V \\
B_0 \> = >> B_0
\end{CD}
$$
such that in a neighborhood of ${g'}^{-1}_{0} \tau^{-1}
(\Sigma_i)$ the following conditions hold true:
\begin{enumerate}
\item[i)] $X'_0$ is non-singular and $\tau'$ induces a
birational morphism $X'_0 \to X_0 \times_{W_0} W'_0$, biregular
over $\tau^{-1} (W_0 - \Sigma^{(0)})$.
\item[ii)] $f'_0 = h'_0 \circ g'_0$ is smooth, outside of a
finite number of points $t_1, \ldots , t_k$.
\item[iii)] For each of the points $t_j$ in ii), there exists a
factorization
$$
\begin{CD}
X'_0 \> \sigma' >> X''_0 \> \delta' >> X_0 \\
\V g'_0 VV \V V g''_0 V \V V g_0 V \\
W'_0 \> \sigma >> W''_0 \> \delta >> W_0 \\
\V h'_0 VV \V V h''_0 V \V V h_0 V \\
B_0 \> = >> B_0 \> = >> B_0
\end{CD}
$$
with $h''_0 \circ g''_0$ smooth in $\sigma' (t_j)$ and $\sigma'$
finite over a neighborhood of $t_j$.
\end{enumerate}
\end{lemma}

\begin{proof}
In what follows we will work locally in $U$, but by abuse of notations
we will write $U = W_0$ and $\Sigma = \Sigma_i$.

Consider first the case where the reduced general fibre of
$\Gamma = g^{-1}_{0} (\Sigma) \to \Sigma$ is a smooth elliptic
curve. Let us assume for a moment that $W''_0 \to W_0$ has
ramification order $m$, the multiplicity of $\Gamma$.
Then the normalization
$$\tilde{g} : \tilde{X}_0 \>>> W''_0 \mbox{ \ \ \ of \ \ \ }
pr_2: X_0 \times_{W_0} W''_0 \>>> W''_0$$
is \'etale over $X_0$, hence smooth over $B_0$. However
$\Gamma_{\red} \simeq \tilde{\Gamma} = \tilde{g}^{-1}
(\Sigma'')$, for $\Sigma'' = \tau^{-1} (\Sigma)_{\red}$ might be
singular. The fibres of $\tilde{\Gamma} \to \Sigma''$ are smooth
elliptic curves or reduced Newton-polygons. Therefore
$\tilde{\Gamma}$ has at most rational Gorenstein singularities.

Let $q$ be one of those singularities, $p = \tilde{g_0} (q)$. We
choose local parameters $(x,y)$ on $W_0$ in $p$, such that
$\Sigma $ is the zero set of $y$ and such that $x$ is the
pullback of a local parameter on $B_0$ in $h(p)$. So the
branched cover $\tau : W''_0 \to W_0$ is locally given by
$\tau^* y = w^m$ and $\tau^* x = x''$, for parameters $w$ and
$x''$ on $W''_0$.

Considering the projection given by $w$ and the morphism induced
in a neighborhood of $q$ in $\tilde{X}_0$, we obtain a family
of surfaces with a smooth general fibre and with an isolated
rational Gorenstein singularity in the special fibre $w
=0$. By \cite{Bri} such a singularity allows a
simultaneous resolution, after taking a further branched
covering, totally ramified over $w =0$. Hence, replacing $m$ by
$m \cdot \nu$ for some $\nu = \nu (q)$ depending on $q$, we may
assume that the normalization $\tilde{X}_0$ of $X_0 \times_{W_0}
W''_0$ has a small resolution $\pi : X''_0 \to \tilde{X}_0$. In
particular, ${g''}^{-1}_{0}  (\Sigma'' ) = \Gamma''$ is smooth and
the fibres of $\Gamma'' \to \Sigma''$ are reduced curves with at
most ordinary double points as singularities, at least over a
neighborhood of the given point.

To do this simultaneously for all points over $\Gamma$, we have
to choose $\mu$ to be divisible by $m$ and by $\nu (q)$ for all
singular points in $\Gamma_{\red}$. For any multiple $M$ of
$\mu$ let $W'_0 \to W_0$ be the corresponding covering. Locally,
for the point $q$ considered
above, $X'_0$ can be chosen to be the covering of $X''_0$,
totally ramified along $\Gamma''$ of order $\frac{M}{m \cdot \nu
(q)}$. Since $\Gamma''$ is non-singular, $X'_0$ is non-singular,
and i) holds true. The conditions ii) and iii) follow from the
construction of $X'_0$.

If the general fibre of $\Gamma_{\red} \to \Sigma$ is a Newton
polygons of length $b$ (hence of type $I_b$), the construction
of $X'_0$ is quite similar. All fibres of $\Gamma_{\red} \to
\Sigma$ are of type $I_a$, for $a \geq b$. Again we start with
$W''_0 \to W_0$, totally ramified of order $m$ and with the
normalization $\tilde{g} : \tilde{X}_{0} \to W''_0$. The
non-smooth locus of $\Gamma_{\red} = \tilde{\Gamma} \to
\Sigma''$ consists of $b$ disjoint sections, say $L_1, \ldots ,
L_b$ and of a finite number of points in $\tilde{\Gamma} - (L_1
\cup \ldots \cup L_b)$. For the latter, the argument given above
works. In fact, if $q$ is one of the isolated points,
$\tilde{\Gamma} - ( L_1 \cup \ldots \cup L_b)$ has again a
rational Gorenstein singularity in $q$, and choosing a larger
covering there exists a simultaneous resolution.

Along $L_i$, the morphism $\tilde{\Gamma} \to \Sigma''$ is
equi-singular. Hence for $W'_0 \to W''_0$ totally ramified over
$\Sigma''$, we can simultaneously resolve the singularities in
the normalization of $X''_0 \times_{W'_0} W''_0$ which are lying
over $L_i$.
\end{proof}

The morphism $f'_0$ constructed in \ref{multiple-locus} is
smooth outside of a finite number of points $t_1, \ldots , t_k$.
As in \ref{defect}, for $T= \{ t_1, \ldots , t_k\}$, let
$(\Omega^{2}_{X'_0/B_0})' = \omega'_{X'_0/B_0} $ be the image of
$\Omega^{2}_{X'_0} \to \omega_{X'_0/B_0}.$

\begin{corollary} \label{prime}
For $M$ and $\tau'$ as in \ref{multiple-locus} and for 
$U'=(g_0\circ \tau')^{-1}(U)$, one has a natural inclusion
$$
\varphi : {\tau'}^{*} \omega_{X_0/B_0}|_{U'} \>>> \omega'_{X'_0/B_0}|_{U'}.
$$
\end{corollary}

\begin{proof}
Both sheaves are subsheaves of $\omega'_{X'_0/B_0}$, hence in
order to show that the inclusion ${\tau'}^{*} \omega_{X_0/B_0} \to
\omega_{X'_0/B_0}$ factors through $\omega'_{X'_0/B_0}$, we can
argue locally in a neighborhood of $t_j \in T$.

Using the notation from \ref{multiple-locus} iii), one has
(locally in a neighborhood of $t_j$) a diagram
$$
\begin{CD}
{\sigma'}^{*} \Omega^{2}_{X''_0} \,\>>> \hspace{-.55cm}
\to\hspace{.15cm} {\sigma'}^{*}
\omega_{X''_0/B_0 } \< \subset << {\delta'}^{*} {\sigma'}^{*}
\omega_{X_0/B_0} \\
\V V \subset V \V V \subset V \noarr \\
\Omega^{2}_{X'_0} \,\>>> \hspace{-1.1cm}
\to\hspace{.15cm} \omega'_{X'_0/B_0} \> \subset >>
\omega_{X'_0/B_0}
\end{CD}
$$
and \ref{prime} holds true.
\end{proof}

\begin{remark} \label{prime2}
The corollary \ref{prime} obviously remains true in the
following situation: Let $\Theta \subset W_0$ be a section, such
that $g^{-1}_{0} (\Theta)$ is non-singular and such that
$g^{-1}_{0} (\Theta) \to \Theta $ is smooth outside of a finite number
of points. Let $W'_0 \to W_0$ be (locally near $\Theta$) a
covering, totally ramified of order $M$ along $\Theta$ and unramified
elsewhere. Then $X'_0 = X_0 \times_{W_0} W'_0$ is non singular,
$X'_0 \to B_0$ is smooth outside of a finite number of points,
and there is (locally) a natural inclusion
$$
\varphi : pr^{*}_{1} \omega_{X_0/B_0} \>>> \omega'_{X'_0/B_0} .
$$
\end{remark}
\section{The proof of \ref{main-theorem} for families of
elliptic surfaces}

Let $f: X \to B$ be a family of minimal elliptic surfaces,
smooth over $B_0 = B-S$ and with $\kappa (X_b) =1$ for $b \in B_0$.
At the beginning of section~4 we proved that
$f$ is birationally isotrivial, in case $B=B_0$ is an elliptic curve.
Hence we will restrict ourselves to the case $B_0=\C^*$, in
this section.

By \ref{easy} we only have to consider families of elliptic
surfaces with one or two multiple
fibres. Nevertheless, since the arguments used here apply to all
other cases as well, we will not make this restriction.

For the relative Iitaka fibration $X_0 \> g_0 >> W_0 \> h_0 >> B_0$,
constructed in section~3, $h_0: W_0 \to B_0$ is a smooth
family of curves. By \ref{trivial} we may write $W_0 = C
\times B_0$, and \ref{multiple} allows to assume that the
multiple locus $\Sigma^{(0)}$ in $W_0$ consists of disjoint
sections. Replacing $C$ by an
\'etale cover, and using \ref{isotrivial-cover} we are allowed
to assume that $g_b : X_b \to W_b=C$ has more than 2 multiple
fibres, provided $g (W_b) \geq 1$.

The further construction depends on the number and type of the singular
fibres: \\[.3cm]
\underline{Case I:} The degree $r$ of $\Sigma^{(0)}$ over $B_0$
is larger than or equal to 2. \\[.2cm]
\underline{Case II:} If $r =1$, we have to find a second
section $\Theta$, with $g^{-1}_{0} (\Theta) \to \Theta$ smooth
outside of a finite number of points, and with $g^{-1}_{0} (\Theta)$
non-singular.

For $g(W_b)> 0$, we were allowed
to assume that $r>1$, hence it is sufficient to construct $\Theta$
for $W_0 \simeq \P^1 \times B_0$. Then the canonical
bundle formula and the assumption $\kappa (X_b) =1$ imply that
$\deg ({g_b}_*\omega_{X_b/W_b}) = \chi (\sO_{X_b} ) \geq 2.$
Depending on the singular fibres of $X_b \to W_b$, for $b\in
B_0$ in general position, we distinguish two subcases: \\[.1cm]
\underline{Case II a:} All non-multiple singular fibres of $X_b
\to W_b$ are reduced. Hence the only singular fibres are of type
${}_mI_n$, $II, III$ or $IV$. We choose for $\Theta$ a
general section of $W_0 \to B_0$, not meeting $\Sigma_1$.\\[.1cm]
\underline{Case II b:} If $X_b \to W_b$ has singular fibre of
type $I^{*}_{n}$, $II^*$, $III^*$ or $IV^*$, we have to be more 
careful. By \ref{weyl} the $I^{*}_{n}$, $II^*$, $III^*$ and
$IV^*$-loci are disjoint in $W_0$, and \'etale over $B_0$.
Replacing $B_0$ again by an \'etale cover we find sections 
$D_{s+1}, \ldots , D_{\ell'}$, corresponding to singular fibres
of type $I^{*}_{n}$, $II^*$, $III^*$ or $IV^*$.
The isomorphism $W_0 \simeq \P^1\times B_0$ can be chosen, such
that $pr_1(\Sigma_1)$, and $pr_1(D_i)$ are points in $\P^1$,
necessarily distinct. 

In fact this is obvious, for $\ell'=s+1$ or $\ell'=s+2$.
If $\ell' > s+2$, we choose the isomorphism such that
$pr_1(\Sigma_1)$, $pr_1(D_{s+1})$ and $pr_1(D_{s+2})$ are points.
Since $D_i \simeq \C^*$, for $i=s+1, \ldots , \ell'$, and since
the $D_i$ can not meet each other or $\Sigma_1$, the restriction
$pr_1|_{D_i}$ can not be dominant, for $i>s+2$. 

For $\Theta$ we choose the fibre of $pr_1$ over a point in
general position in $\P^1$.\\[.3cm] 
Let $\mu_i$ be the number, attached to the component $\Sigma_i$
in \ref{multiple-locus} and let $M= {\rm lcm} \{ \mu_1, \ldots ,\mu_r
\}$. We choose a covering $\tau_0 : W'_0 \to W_0$,
totally ramified of order $M$ over $\Sigma_1 + \Theta$, in case
II, or ramified in each component of
$\displaystyle \sum^{r}_{i=1} \tau_0^* \Sigma_i$ of order $M$,
in case I, and we assume $\tau_0$ to be unramified elsewhere.

Such coverings exist by \cite{FK}, IV.9.12, for example. As
explained in \cite{EV}, 3.5 and 3.15, 
they can also be obtained by taking the $m$-root out of divisors
$A$, with $\sO_{W_0}(A)$ the $m$-th power of an invertible sheaf.  

In case II, choose $A= \Sigma_1 + (M-1) \cdot \Theta$, and in
case I, if $r$ is even,
$$
A=\sum_{i=1}^{r/2} \Sigma_{2i-1} + (M-1)\cdot \Sigma_{2i}
$$
will give the covering needed.
If $r$ and $M$ are both odd, one can take the $M$-th root
out of 
$$
A= \Sigma_{r-2} + \Sigma_{r-1} + (M-2)\cdot \Sigma_r +
\sum_{i=1}^{(r-3)/2} \Sigma_{2i-1} + (M-1)\cdot \Sigma_{2i}.
$$
If $M$ is even, and $r$ odd, take first the $M$-th root out
of $A=\Sigma_1 + (M-1)\cdot \Sigma_2$. On the covering obtained there are 
$M\cdot (r-2)$ points left, and we proceed as in the first step.

Let $X'_0 \to W'_0$ be the model from
\ref{multiple-locus} over a neighborhood of $\tau^{-1}_{0}
(\Sigma^{(0)})$, in both cases, and equal to $X_0
\times_{W_0} W'_0$ over $\Theta$, in case II.
We constructed a non-singular variety $X'_0$ and projective
morphisms
$$
\begin{CD}
X'_0 \> \tau'_0 >> X_0 \\
\V g'_0 VV \V V g_0 V \\
W'_0 \> \tau_0 >> W_0
\end{CD}
$$
such that $f'_0 = h_0 \circ \tau_0 \circ g'_0$ is smooth outside
of a finite number of points, and such that locally in those
points \ref{prime} (see also \ref{prime2}) holds. This
remains true if we replace $B_0$ by further \'etale coverings.
Hence we may choose non-singular projective compactifications
$$
\begin{CD}
X' \> \tau' >> X \\
\V g' VV \V V g V\\
W' \> \tau >> W \\
\V h' VV \V V h V \\
B \> = >> B
\end{CD}
$$
such that $f' = h' \circ g'$ as well as $f = h \circ g$ satisfy
the conditions stated in \ref{prepared} (for $B'=B$).

{\it Assume that $f : X \to B$ is not birationally isotrivial.}
Then \ref{positivity} implies that $\kappa(\omega_{X/B}) = 2$. Using
the notations from \ref{prepared}, viii), (with all the ${}'$ omitted), one
has
$$
f_* \omega^{\nu}_{X/B} = h_* \Bigl(\omega^{\nu}_{W/B} \otimes \sO_W
\Bigl( \sum^{r}_{i=1} \Bigl[ \frac{\nu \cdot (m_i-1)}{m_i}
\Bigr] \cdot \Sigma_i \Bigr) \otimes \delta^{\nu} \Bigr).
$$
Hence, if $\nu$ is a multiple of ${\rm lcm}\{ m_1, \ldots , m_r\}$, the sheaf
$$
\sL^{(\nu)} = \omega^{\nu}_{W/B} \otimes \sO_W \Bigl(
\sum^{r}_{i=1} \frac{\nu \cdot (m_i-1)}{m_i} \Sigma_i \Bigr)
\otimes \delta^{\nu}
$$
will be ample with respect to $W_0$, and, for some
effective divisor $R$ on $X$ supported in $f^{-1} (S)$, one has
$$
g^* \sL^{(\nu)} = \omega_{X/B} (-R)^{\nu} .
$$
Since $\tau$ is totally ramified over $\Sigma_i$ of order
divisible by $m_i$, there is an invertible sheaf $\sL'$ on $W'$,
ample with respect to $W'_0$, with $\tau^* \sL^{(\nu)} =
{\sL'}^{\nu}$. Moreover ${g'}^{*} \sL' = {\tau'}^{*} \omega_{X/B}
(-R)$, and $\sL'$ is a subsheaf of $g'_* \omega_{X'/B}$.

Let $\gamma': J' \to W'$ be a Jacobian fibration, as considered in
section~4. By \ref{compare} one has
$$
\gamma'_* \omega^{\nu}_{J'/W'} = g'_* \omega^{\nu}_{X'/W'},
$$
for all $\nu \geq 1$. Replacing a last time $B_0$ by an \'etale
cover, \ref{prepared2} allows to assume that $J' \to B$
satisfies the conditions stated in \ref{prepared} (over $B=B'$).

Since $\gamma'_0 : J'_0 \to W'_0$ is locally in the \'etale
topology isomorphic to $g'_0 : X'_0 \to W'_0$, the morphism
$J'_0 \to B_0$ is again smooth outside of a finite subset $T$.

$J'_0 \to W'_0$ has a zero-secton with image $\Pi_0$. Writing $\psi'=h'
\circ \gamma'$ and ${\psi'}^{-1} (S) = \nabla'$, we may assume that the
closure $\Pi$ of $\Pi_0$ is non-singular and that $\nabla' +
\Pi$ is a normal crossing divisor.

As in \ref{defect} one defines
\begin{gather*}
\Omega^{2}_{X'/B} (\log\Delta')' = \im (\Omega^{2}_{X'}
(\log\Delta') \>>> \Omega^{2}_{X'/B} (\log \Delta')^{\sim} ),
\mbox{ \ \ \ and}\\
\Omega^{2}_{J'/B} (\log\nabla')' = \im (\Omega^{2}_{J'}
(\log\nabla') \>>> \Omega^{2}_{J'/B} (\log \nabla')^{\sim} ).
\end{gather*}
$\omega'_{X'/B}$ denotes the subsheaf of $\omega_{X'/B}$, generated by
$\Omega^{2}_{X'/B}  (\log\Delta')'$ and by $\omega_{X'/B}$,
restricted to a neighborhood of $\Delta'$. Correspondingly
$\omega'_{J'/B}$ is generated by the sheaves $\Omega^{2}_{J'/B}
(\log\nabla')'$ and $\omega_{J'/B}|_{J'-T}$.

Since $\gamma'_0 : J'_0 \to W'_0$ and $g'_0 : X'_0 \to W'_0$ are
locally isomorphic in the \'etale topology, the natural
isomorphism
$$
\gamma'_{0*} \omega_{J'/B} \> \simeq >> g'_{0*} \omega_{X'/B}
$$
induces an isomorphism
\begin{equation}\label{iso1}
\gamma'_{0*} \omega'_{J'/B} \> \simeq >> g'_{0*} \omega'_{X'/B}.
\end{equation}
By \ref{prime} (see also \ref{prime2}),
${g'}^{*} \sL' = {\tau'}^{*} \omega_{X/B} (-R)$
is contained in $\omega'_{X'/B}$. Hence $\sL'$ lies in
$g'_* \omega'_{X'/B}$, and, using the isomorphism (\ref{iso1})
one finds an injection
$$
{\gamma'}^*{\sL'} \> \subset >> \omega'_{J'/B}.
$$
The second sheaf contains $\Omega^2_{J'/B}(\log \nabla')'$ and
for some effective divisor $\Upsilon$ on $J'$, supported in $\nabla'$,
\begin{equation}\label{Ups}
\Omega^2_{J'/B}(\log \nabla')'\cap
{\gamma'}^*\sL'={\gamma'}^*\sL\otimes \sO_{J'}(-\Upsilon).
\end{equation}
By the last condition in \ref{prepared}, for all $\nu>0$,
$$
\psi'_*\Omega^2_{J'/B}(\log \nabla')^\nu =
\psi'_*\omega_{J'/B}^\nu,
$$
hence
$$
\psi'_*({\gamma'}^*{\sL'}^\nu \otimes \sO_{J'}(-\nu\cdot
\Upsilon)) = h'_*{\sL'}^\nu.
$$
Altogether $\sL'$, $\Pi$ and $\Upsilon$ satisfy the assumptions
made in \ref{van-ass}, for $J'\to W' \to B$. By \ref{van}
$$
H^0 (J', \Omega^{2}_{J'/B} (\log \nabla')' \otimes
{\gamma'}^* \sL^{-1} \otimes \sO_{J'}(\Upsilon)) = 0,
$$
contradicting the choice of $\Upsilon$ in (\ref{Ups}).
\qed
\bibliographystyle{plain}

\end{document}